\documentclass[12pt]{article}
\usepackage{amsmath}
\usepackage{amssymb,amsfonts}
\usepackage{amscd}

\newtheorem{theorem}{Theorem}[section]
\newtheorem{lemma}[theorem]{Lemma}
\newtheorem{proposition}[theorem]{Proposition}

\newtheorem{example}[theorem]{Example}
\newtheorem{definition}[theorem]{Definition}
\newtheorem{corollary}[theorem]{Corollary}
\newtheorem{remark}[theorem]{Remark}
\newenvironment{proof}{\bf Proof. \rm}{$\Box$}
\newcommand{\be}{\begin{equation}}
\newcommand{\ee}{\end{equation}}

\begin{document}

\title{Quantum Markov Semigroups\\
(Product Systems and Subordination)}
\author{Paul S. Muhly\thanks{%
Supported in part by grants from the National Science Foundation and from
the U.S.-Israel Binational Science Foundation.} \\
Department of Mathematics\\
University of Iowa\\
Iowa City, IA 52242\\
e-mail: muhly@math.uiowa.edu \and Baruch Solel\thanks{%
Supported in part by the U.S.-Israel Binational Science Foundation and by
the Fund for the Promotion of Research at the Technion.} \\
Department of Mathematics\\
Technion\\
32000 Haifa, Israel\\
e-mail: mabaruch@techunix.technion.ac.il}
\date{}
\maketitle

\begin{abstract}
We show that if a product system comes from a quantum Markov semigroup, then
it carries a natural Borel structure with respect to which the semigroup may
be realized in terms of a measurable representation. We show, too, that the
dual product system of a Borel product system also carries a natural Borel
structure. We apply our analysis to study the order interval consisting of
all quantum Markov semigroups that are subordinate to a given one.
\end{abstract}

\section{Introduction\label{Intro}}

A \emph{quantum Markov semigroup} is a semigroup $\{\Theta _{t}\}_{t\geq 0}$
of completely positive, normal linear maps on a von Neumann algebra $M$ such
that $\Theta _{0}$ is the identity mapping on $M$ and such that the map $%
t\rightarrow \Theta _{t}(a)$ from $[0,\infty )$ to $M$ is continuous with
respect to the ultraweak topology on $M$ for each $a\in M$. If each $\Theta
_{t}$ is a unital map, we shall say that the semigroup is \emph{unital.} In
\cite{MSQMP}, we showed that if $M$ is countably decomposable, then it is
possible to dilate a unital quantum Markov semigroup to an \textquotedblleft
endomorphic\textquotedblright\ semigroup in the following sense. Represent $%
M $ faithfully on a separable Hilbert space, say $H$. Then one may find:
another separable Hilbert space $K$, an isometric embedding $u_{0}$ of $H$
into $K$, a von Neumann algebra $\mathcal{R}$ in $B(K)$ and a semigroup of
unital \emph{endomorphisms }of $\mathcal{R}$, $\{\alpha _{t}\}_{t\geq 0}$,
such that $u_{0}Mu_{0}^{\ast }$ is a full corner in $\mathcal{R}$, meaning
that the central support of $u_{0}u_{0}^{\ast }$ in $\mathcal{R}$ is $I$; $%
\{\alpha _{t}\}_{t\geq 0}$ is a quantum Markov semigroup; and such that the
two (equivalent) equations are satisfied for all $T\in M$, all $S\in
\mathcal{R}$, and all $t\geq 0$:
\begin{equation*}
\Theta _{t}(T)=u_{0}^{\ast }\alpha _{t}(u_{0}Tu_{0}^{\ast })u_{0}
\end{equation*}%
and
\begin{equation*}
\Theta _{t}(u_{0}^{\ast }Su_{0})=u_{0}^{\ast }\alpha _{t}(S)u_{0}\text{.}
\end{equation*}%
Semigroups of unital endomorphisms of a von Neumann algebra, such as $%
\{\alpha _{t}\}_{t\geq 0}$, which are also quantum Markov semigroups are
known in the literature as $E_{0}$-semigroups and were first defined and
investigated by Powers \cite{rP88} and Arveson \cite{Arv89}.

The proof of our dilation theorem proceeded by expressing $\{\Theta
_{t}\}_{t\geq 0}$ in terms of a representation of a product system $%
\{E(t)\}_{t\geq 0}$ of $W^{\ast }$-correspondences over the commutant of $M$%
, $M^{\prime }$. (Definitions and further details will be given below.) In
\cite{MSQMP} we attended only to the algebraic structure of $\{E(t)\}_{t\geq
0}$; that is all that we needed there. In this sequel, our primary objective
is to show how to put a Borel structure on $\{E(t)\}_{t\geq 0}$ and to
relate the Borel structure to continuity properties of its representations
(see Theorems \ref{isomprod}, \ref{meas}, and \ref{cont}). Our approach to
dilating quantum Markov semigroups is closely related to the approach taken
by Bhat and Skeide \cite{BS00}. Indeed, the two approaches are
\textquotedblleft dual\textquotedblright\ in a sense made precise in
Skeide's survey \cite{mS03}. Our analysis shows that the product system that
Bhat and Skeide construct is also Borel. (See Theorem \ref{dualprodsys}.)

As an application of our analysis, we study quantum Markov semigroups that
are \textquotedblleft subordinate\textquotedblright\ to a given one. If $%
\{\Theta _{t}\}_{t\geq 0}$ and $\{\Psi _{t}\}_{t\geq 0}$ are two quantum
Markov semigroups on the same von Neumann algebra, then we say that $\{\Psi
_{t}\}_{t\geq 0}$\ is subordinate to $\{\Theta _{t}\}_{t\geq 0}$, if $\Theta
_{t}-\Psi _{t}$ is completely positive for all $t\geq 0$. We show that the
subordinates of a given semigroup $\{\Theta _{t}\}_{t\geq 0}$ depend only on
the product system $\{E(t)\}_{t\geq 0}$ associated to the semigroup and not
directly on the semigroup itself, provided it comes from a \textquotedblleft
injective representation\textquotedblright\ of $\{E(t)\}_{t\geq 0}$ in a
sense that we define below (see Theorem \ref{subordinatetheorem}). Thus, if $%
\{\Theta _{t}\}_{t\geq 0}$ and $\{\Psi _{t}\}_{t\geq 0}$ both come from
injective representations of $\{E(t)\}_{t\geq 0}$, then their order
intervals\ of subordinates are order isomorphic. These isomorphism results
were proved by Bhat \cite[Section 5]{bB01} and Powers \cite[Theorems 3.4 and
3.5]{Po03} in the case when $M=B(H)$ by entirely different means and made no
use of product systems. Another result similar to ours was proved by Bhat
and Skeide in \cite[Theorem 14.3]{BS00}. Indeed, from a purely algebraic
perspective, they come to the same conclusion, but from a perspective that
is dual to ours, in a sense that we shall discuss at several points below.
Our contribution is to deal with continuity properties of the Markov
semigroups and the Borel structures on the associated product systems and to
show in terms of the representation theory for product systems that we
developed in \cite{MSQMP} that the order interval of subordinates to a given
semigroup is an artifact of the product system of which it is a
representation and does not require the dilation. More explicitly, we define
the notion of a (positive contractive) cocycle for a product system
(Definition \ref{reducedcocycles}) and show in Theorem \ref%
{subordinatetheorem} that if $\{\Theta _{t}\}_{t\geq 0}$ comes from an
injective representation of a product system $\{E(t)\}_{t\geq 0}$, then
every quantum Markov semigroup that is subordinate to $\{\Theta
_{t}\}_{t\geq 0}$ is given by a cocycle for $\{E(t)\}_{t\geq 0}$. In the
works of Bhat and Powers just cited, the approach is to pass from $\{\Theta
_{t}\}_{t\geq 0}$ to its minimal endomorphic dilation $\alpha $ and to
express subordinates of $\{\Theta _{t}\}_{t\geq 0}$ in terms of so-called
local cocycles for $\alpha $. As we shall show in Proposition \ref%
{cocycleequal}, if $\alpha $ comes from an isometric representation of a
product system $\{E(t)\}_{t\geq 0}$, then there is a natural bijection
between cocycles for $\{E(t)\}_{t\geq 0}$ and local cocycles for $\alpha $.
We want to emphasize, however, that given $\{\Theta _{t}\}_{t\geq 0}$ it is
not necessary to pass to its minimal endomorphic dilation $\alpha $ in order
to analyze the subordinates and, in fact, our analysis works without the
assumption that $\{\Theta _{t}\}_{t\geq 0}$ has such a dilation. In
particular, it allows for a direct comparison of the subordinates of two
different semigroups (Corollary \ref{orderisom}).\medskip
\newpage
\noindent\textbf{Acknowledgment }

We thank the referee for some very stimulating comments that helped us to
improve the first draft of this paper substantially. In particular, some
bothersome loose ends have been cleared up and the arguments are a lot
clearer.\bigskip

\noindent \textbf{Conventions and Notation}

The purely algebraic aspects of our analysis require no separability
hypotheses. However, when we deal with continuity properties of semigroups
or Borel structures on product systems, we will need to assume that our von
Neumann algebras are countably decomposable, meaning that they can be
faithfully represented on a separable Hilbert space, and we will need to
assume that representation Hilbert spaces, when they arise, are separable.
This will guarantee, among other things that all represented von Neumann
algebras have countably decomposable commutants.

If $S$ is a nonempty subset of a Hilbert space $H$, then we shall denote the
closed linear span of $S$ by $[S]$.

\section{Completely Positive Maps,\newline
$W^{\ast }$-correspondences, and their Representations\label%
{CPMapsWsCorrReps}}

We collect in this section facts about $W^{\ast }$-correspondences and their
representations that we need in order to study a normal contractive
completely positive map on a von \ Neumann algebra. Most of the material we
discuss may be found in \cite{MSQMP}. However, we need some refinements of
the theory presented there and we want to highlight certain features of it.
In a bit more detail, we show how a completely positive (contractive normal)
map on a von Neumann algebra gives rise to two $W^{\ast }$-correspondences,
the principal one for us being the Arveson-Stinespring correspondence, which
is a correspondence over the commutant of the von Neumann algebra. The
completely positive map also gives a representation of the correspondence.
Conversely, each representation of a correspondence gives rise to a
completely positive map. Our primary objective, from the technical point of
view, is to analyze when the transformation from a completely positive map
to correspondence and representation is reversable.

The theory of $W^{\ast }$-correspondences is based on the theory of Hilbert $%
C^{\ast }$-modules. We shall follow Lance \cite{L94} and Paschke \cite{Pa73}
for notation and the parts of that theory that we shall use. Let $A$ be a $%
C^{\ast }$-algebra and let $E$ be a right module over $A$ endowed with a
bi-additive map $\langle \cdot ,\cdot \rangle :E\times E\rightarrow A$
(referred to as an $A$-valued inner product) such that, for $\xi ,\eta \in E$
and $a\in A$, $\langle \xi ,\eta a\rangle =\langle \xi ,\eta \rangle a$, $%
\langle \xi ,\eta \rangle ^{\ast }=\langle \eta ,\xi \rangle $, and $\langle
\xi ,\xi \rangle \geq 0$, with $\langle \xi ,\xi \rangle =0$ only when $\xi
=0$. Also, $E$ is assumed to be complete in the norm $\Vert \xi \Vert
:=\Vert \langle \xi ,\xi \rangle \Vert ^{1/2}$. We write $\mathcal{L}(E)$
for the space of adjointable, and therefore continuous $A$-module maps on $E$%
. It is known to be a $C^{\ast }$-algebra. The $C^{\ast }$-module $E$ is
said to be \emph{self-dual }in case every adjointable $A$-module map from $E$
to $A$ is given by an inner product with an element of $E$. If $A$ is a von
Neumann algebra and $E$ is self-dual, then every continuous module map is
adjointable and $\mathcal{L}(E)$ is a von Neumann algebra. (See \cite[%
Corollary 3.5 and Proposition 3.10]{Pa73}.) A $C^{\ast }$-\emph{%
correspondence} over a $C^{\ast }$-algebra $A$ is a Hilbert $C^{\ast }$%
-module $E$ over $A$ that is endowed with a structure of a left module over $%
A$ via a nondegenerate $\ast $-homomorphism $\varphi :A\rightarrow \mathcal{L%
}(E)$. When dealing with a specific $C^{\ast }$-correspondence, $E$, over a $%
C^{\ast }$-algebra $A$, it will be convenient to suppress the $\varphi $ in
formulas involving the left action and simply write $a\xi $ or $a\cdot \xi $
for $\varphi (a)\xi $. \ This should cause no confusion in context.

\begin{definition}
Let $M$ be a von Neumann algebra and let $E$ be a Hilbert $C^{\ast }$-module
over $M$. Then $E$ is called a \emph{Hilbert }$W^{\ast }$\emph{-module} over
$M$ in case $E$ is self-dual. The module $E$ is called a $W^{\ast }$\emph{%
-correspondence over} $M$ in case $E$ is a self-dual $C^{\ast }$%
-correspondence over $M$ such that the $\ast $-homomorphism $\varphi
:M\rightarrow \mathcal{L}(E)$ giving the left module structure on $E$ is
normal.
\end{definition}

As we mentioned above, in the portions of this paper that deal with
measurable product systems we will restrict our attention to countably
decomposable von Neumann algebras and normal representations of them on
separable Hilbert spaces. No separability assumptions are necessary in this
section.

One of the most important examples of $W^{\ast }$-correspondences arises
from completely positive normal maps on a von Neumann algebra, as in the
following example.

\begin{example}
\label{cpmod}Suppose that $\Theta $ is a normal, contractive, completely
positive map on a von Neumann algebra $M$. Then we can associate with it the
correspondence $M\otimes _{\Theta }M$ obtained by defining on the algebraic
tensor product $M\otimes M$, the $M$-valued inner product $\langle a\otimes
b,c\otimes d\rangle =b^{\ast }\Theta (a^{\ast }c)d$ and taking the selfdual
completion (\cite[Theorem 3.2]{Pa73}). The bimodule structure is defined by
left and right multiplications. This correspondence was used by Popa \cite%
{P86}, Mingo \cite{Mi89}, Anantharaman-Delaroche \cite{AnD90} and others to
study the map $\Theta $. Bhat and Skeide refer to $M\otimes _{\Theta }M$ as
the GNS-module determined by $\Theta $ in \cite{BS00}; we shall call it the
\emph{GNS-correspondence} determined by $\Theta $. The reason for the
terminology is evident when one notes that there is a preferred vector $\xi $
in $M\otimes _{\Theta }M$, which is the image of $1\otimes 1$. The
completely positive map $\Theta $ may be expressed in terms of $\xi $ via
the formula $\Theta (a)=$ $\langle \xi ,a\xi \rangle $, $a\in M$, which is
analogous to the GNS representation determined by a state on a $C^{\ast }$%
-algebra. If $\Theta $ is an endomorphism, we can use the map $a\otimes
_{\Theta }b\mapsto \Theta (a)b$ to identify this correspondence with the
correspondence $_{\Theta }M$ which is defined to be the subspace $\Theta
(I)M $ of $M$ with right action of $M$ given by multiplication (on the
right), left action by $\varphi (a)=\Theta (a)$ and the inner product
induced from multiplication on $M$, $\langle c,d\rangle =c^{\ast }d$.
\end{example}

Given two $W^{\ast }$-correspondences $E$ and $F$ over $M$, then the
balanced tensor product carries a natural inner product, which is defined by
the formula
\begin{equation*}
\langle \xi _{1}\otimes \zeta _{1},\xi _{2}\otimes \zeta _{2}\rangle
:=\langle \zeta _{1},\varphi _{F}(\langle \xi _{1},\xi _{2}\rangle )\zeta
_{2}\rangle
\end{equation*}%
and its completion in the $\sigma $-topology of \cite{BDH88} is a $W^{\ast }$%
-correspondence over $M$, where, for $a,b\in M$ and $\xi \in E$ and $\zeta
\in F$, $\varphi _{E\otimes F}(a)(\xi \otimes \zeta )b=\varphi _{E}(a)\xi
\otimes \zeta b$.

\begin{definition}
\label{repres}Let $E$ be a $W^{\ast }$-correspondence over a von Neumann
algebra $N$. Then:

\begin{enumerate}
\item A \emph{completely contractive covariant representation }of $E$ on a
Hilbert space $H$ is a pair $(T,\sigma)$, where

\begin{enumerate}
\item $\sigma$ is a normal $\ast$-representation of $N$ in $B(H)$.

\item $T$ is a linear, completely contractive map from $E$ to $B(H)$ that is
a bimodule map in the sense that $T(a\xi b)=\sigma (a)T(\xi )\sigma (b)$, $%
\xi \in E$, and $a,b\in N$.
\end{enumerate}

\item A completely contractive covariant representation $(T,\sigma)$ of $E$
in $B(H)$ is called \emph{isometric }in case
\begin{equation}  \label{isometric}
T(\xi)^{\ast}T(\eta)=\sigma(\langle\xi,\eta\rangle)
\end{equation}
for all $\xi,\eta\in E$.
\end{enumerate}
\end{definition}

It should be noted that the operator space structure of $E$ which Definition %
\ref{repres} refers to is that which $E$ inherits when viewed as a subspace
of its linking algebra. We note, too, that Lemma 2.16 of \cite{MSQMP} shows
that if $(T,\sigma )$ is completely contractive covariant representation of $%
E$ on a Hilbert space $H$, then $T$ is continuous in the $\sigma $-topology
of \cite{BDH88} on $E$ and the ultraweak topology on $B(H)$.

Given a $W^{\ast }$-correspondence $E$ over $M$ and a normal representation $%
\sigma $ of $M$ on $H$, we write $E\otimes _{\sigma }H$ for the Hilbert
space obtained as the Hausdorff completion of $E\otimes H$ with respect to
the positive semi-definite sesquilinear form defined by the formula $\langle
\xi \otimes h,\zeta \otimes k\rangle =\langle h,\sigma (\langle \xi ,\zeta
\rangle )k\rangle $. Note that given $S\in \mathcal{L}(E)$ and $R\in \sigma
(M)^{\prime }$, the operator $S\otimes R$, defined by sending $\xi \otimes h$
to $S\xi \otimes Rh$, is a well defined, bounded operator on $E\otimes
_{\sigma }H$. It is easy to see that one obtains ultraweakly continuous $%
\ast $-representations of $\mathcal{L}(E)$ and $\sigma (M)^{\prime }$ on $%
E\otimes _{\sigma }H$ defined by the formulae $S\rightarrow S\otimes I_{H}$
and $R\rightarrow I_{E}\otimes R$. The representation of $\mathcal{L}(E)$, $%
S\rightarrow S\otimes I_{H}$, is called the representation of $\mathcal{L}%
(E) $ \emph{induced }by $\sigma $ and is sometimes denoted $\sigma ^{E}$.
Although it is not standard to do so, we shall call the representation of $%
\sigma (M)^{\prime }$, $R\rightarrow I_{E}\otimes R$, the representation of $%
\sigma (M)^{\prime }$ that is \emph{produced by }$E$. We denote the ranges
of these representations by $\mathcal{L}(E)\otimes I_{H}$ and $I_{E}\otimes
\sigma (M)^{\prime }$, respectively. We then record for the sake of
reference the following lemma which is a restatement of \cite[Theorem 6.23]%
{R74}. For a slick proof that uses a von Neumann algebra version of Brown,
Green and Rieffel's linking algebra of an imprimitivity bimodule \cite{BGR77}%
, we recommend Skeide's note \cite[Proposition 2.2]{mSp03}.

\begin{lemma}
\label{comm} The commutant of the induced algebra $\mathcal{L}(E)\otimes
I_{H}$ in $B(E\otimes _{\sigma }H)$ is the produced algebra $I_{E}\otimes
\sigma (M)^{\prime }$.
\end{lemma}

As we showed in \cite[Lemmas 3.4--3.6]{MSTensor} and in \cite{MSHardy}, if a
completely contractive covariant representation, $(T,\sigma )$, of $E$ in $%
B(H)$ is given, then it determines a contraction $\tilde{T}:E\otimes
_{\sigma }H\rightarrow H$ defined by the formula $\tilde{T}(\eta \otimes
h):=T(\eta )h$, $\eta \otimes h\in E\otimes _{\sigma }H$. The operator $%
\tilde{T}$ intertwines the representation $\sigma $ on $H$ and the induced
representation $\sigma ^{E}\circ \varphi :=\varphi (\cdot )\otimes I_{H}$ on
$E\otimes _{\sigma }H$; i.e.
\begin{equation}
\tilde{T}(\varphi (\cdot )\otimes I)=\sigma (\cdot )\tilde{T}.
\label{covariance}
\end{equation}%
In fact we have the following lemma from \cite[Lemma 2.16]{MSHardy}.

\begin{lemma}
\label{CovRep}The map $(T,\sigma)\rightarrow\tilde{T}$ is a bijection
between all completely contractive covariant representations $(T,\sigma)$ of
$E$ on the Hilbert space $H$ and contractive operators $\tilde{T}%
:E\otimes_{\sigma }H\rightarrow H$ that satisfy equation (\ref{covariance}).
Given such a $\tilde{T}$ satisfying this equation, $T$, defined by the
formula $T(\xi)h:=\tilde{T}(\xi\otimes h)$, together with $\sigma$ is a
completely contractive covariant representation of $E$ on $H$. Further, $%
(T,\sigma)$ is isometric if and only if $\tilde{T}$ is an isometry.
\end{lemma}

An important concept that we shall use is that of duality for $W^{\ast }$%
-corres\-pond\-ences. We shall refer mostly to \cite{MSHardy} and follow the
notation and terminology developed there.

\begin{definition}
\label{dualcorespondence}Let $E$ be a $W^{\ast }$-correspondence over $M$.
Let $\sigma :M\rightarrow B(H)$ be a faithful normal representation of the
von Neumann algebra $M$. Then the $\sigma $\emph{-dual} of $E$, denoted $%
E^{\sigma }$, is defined to be
\begin{equation*}
\{\eta \in B(H,E\otimes _{\sigma }H)\mid \eta \sigma (a)=(\varphi (a)\otimes
I)\eta ,\;a\in M\}.
\end{equation*}
\end{definition}

Thus, by virtue of Definition \ref{repres} and Lemma \ref{CovRep}, we see
that the unit ball of $(E^{\sigma })^{\ast }$ may be identified with the
collection of all covariant representations of $E$ (such that the associated
representation of $M$ is $\sigma $). Moreover, $E^{\sigma }$ has the
structure of a $W^{\ast }$-correspondence over the \emph{commutant }of $%
\sigma (M)$, $\sigma (M)^{\prime }$, as described in the following
proposition.

\begin{proposition}
\label{corres} (\cite[Proposition 2.8]{MSHardy}) With respect to the actions
of $\sigma (M)^{\prime }$ and the $\sigma (M)^{\prime }$-valued inner
product defined as follows, $E^{\sigma }$ becomes a $W^{\ast }$%
-correspondence over $\sigma (M)^{\prime }$: For $a,b\in \sigma (M)^{\prime
} $, and $\eta \in E^{\sigma }$, $a\cdot \eta \cdot b:=(I\otimes a)\eta b$,
and for $\eta ,\zeta \in E^{\sigma }$, $\langle \eta ,\zeta \rangle _{\sigma
(M)^{\prime }}:=\eta ^{\ast }\zeta $.
\end{proposition}

It will be convenient to write $\varphi _{\sigma }$ for the left action of $%
\sigma (M)^{\prime }$ on $E^{\sigma }$, i.e., $\varphi _{\sigma
}(a)T=(I_{E}\otimes a)T$, for $a\in \sigma (M)^{\prime }$ and $T\in
E^{\sigma }$.

\begin{example}
\label{cpdual}If $\Theta $ is a contractive, normal, completely positive map
on a von Neumann algebra $M$ and if $E$ is the GNS correspondence $M\otimes
_{\Theta }M$ of Example~\ref{cpmod}, then, for every faithful representation
$\sigma $ of $M$ on a Hilbert space $H$, the $\sigma $-dual of $E$ is the
space of all bounded operators mapping $H$ into the Stinespring space $K$
(associated with $\sigma \circ \Theta $, which maps $M$ to $B(H)$) that
intertwine the representation $\sigma $ on $H$ and the Stinespring dilation
of $\sigma \circ \Theta $. More specifically, recall that $K$ is the
completion of the algebraic tensor product $M\otimes H$ with respect to $M$%
-valued form, $\langle \cdot ,\cdot \rangle $, defined by the formula $%
\langle a\otimes h,b\otimes k\rangle =\langle h,\sigma \circ \Theta (a^{\ast
}b)k\rangle $. The Stinespring dilation of $\sigma \circ \Theta $ is the
\textquotedblleft induced\textquotedblright\ representation $a\rightarrow
a\otimes I$. Thus $E^{\sigma }$ is the space
\begin{equation}
\{X:H\rightarrow M\otimes _{\sigma \circ \Theta }H|X\sigma (a)=(a\otimes
I)X,\;\;a\in M\}.  \label{etheta}
\end{equation}%
Because of the special role this correspondence plays in this note, we write
$E_{\Theta }^{\sigma }$ for it, or simply $E_{\Theta }$ when $\sigma $ is
understood. Following \cite{MSQMP}, we call $E_{\Theta }^{\sigma }$ or $%
E_{\Theta }$ the \emph{Arveson-Stinespring correspondence} of $\Theta $.
Along with the GNS correspondence, the Arveson-Stinespring correspondence
turns out to be of great importance in the study of a completely positive
map $\Theta $. (See \cite{MSQMP} and \cite{MSCurv} in particular.) If $%
M=B(H) $ and $\sigma $ is the identity representation, $E_{\Theta }^{\sigma
} $ is a Hilbert space and was studied by Arveson \cite{Arv89}. Note also
that, if $\Theta $ is an endomorphism, then $E_{\Theta }^{\sigma }$ can be
identified with the space of all operators on $H$ with ranges in $\sigma
\circ \Theta (I)H$ that intertwine $\sigma $ and $\sigma \circ \Theta $.
\end{example}

Observe that if $\sigma $ is a (faithful) normal representation of $M$ on $H$
and if $\iota $ denotes the identity representation of $\sigma (M)^{\prime }$
on $H$, then for any correspondence $E$ over $M$, we can form $E^{\sigma
}\otimes _{\iota }H$ and we can consider the induced representation $\iota
^{E^{\sigma }}\circ \varphi _{\sigma }$ of $\sigma (M)^{\prime }$ and the
produced representation $a\rightarrow I_{E^{\sigma }}\otimes a$, $a\in
\sigma (M)$, of $\sigma (M)$ on $E^{\sigma }\otimes _{\iota }H$. The
following lemma, which is part of \cite[Lemma 3.8]{MSHardy}, shows that
these representations are unitarily equivalent to the representations of $M$
on $E\otimes _{\sigma }H$ with which we started; i.e. $\iota ^{E^{\sigma
}}\circ \varphi _{\sigma }$ is unitarily equivalent to the produced
representation of $\sigma (M)^{\prime }$ and the produced representation of $%
\sigma (M)$, $\sigma (a)\rightarrow I_{E^{\sigma }}\otimes \sigma (a)$, is
unitarily equivalent to the induced representation $\sigma ^{E}\circ \varphi
$ of $M$.

\begin{lemma}
\label{u} Let $E$ be a $W^{\ast }$-correspondence over a von Neumann algebra
$M$. Let $\sigma $ be a faithful normal representation of $M$ on $H$ and let
$E^{\sigma }$ be the $\sigma $-dual of $E$. Then the map which sends $%
X\otimes h\in E^{\sigma }\otimes _{\iota }H$ to $Xh\in E\otimes _{\sigma }H$
extends to a Hilbert space isomorphism $u:E^{\sigma }\otimes _{\iota
}H\rightarrow E\otimes _{\sigma }H$. Moreover, for $b\in \sigma (M)^{\prime
} $ and $a\in M$,
\begin{equation*}
u(\varphi _{\sigma }(b)\otimes I_{H})=(I_{E}\otimes b)u
\end{equation*}%
and
\begin{equation*}
u(I_{E^{\sigma }}\otimes \sigma (a))=(\varphi _{E}(a)\otimes I_{H})u.
\end{equation*}

In particular, when $E$ is the GNS correspondence associated with a
completely positive map $\Theta $ on $M$, $M\otimes _{\Theta }M$ , $u$ maps $%
E_{\Theta }^{\sigma }\otimes H$ onto $M\otimes _{\Theta }H=E\otimes _{\sigma
}H$.
\end{lemma}

\begin{proof}
The fact that $u$ is a well defined unitary operator can be found in \cite[%
Lemma 3.8]{MSQMP}. The rest follows from the following straightforward
computation. We have, for $b\in \sigma (M)^{\prime },X\in E^{\sigma },h\in H$
and $a\in M$, $u(b\cdot X\otimes h)=b\cdot Xh=(I_{E}\otimes
b)Xh=(I_{E}\otimes b)u(X\otimes h)$ and $u(X\otimes \sigma (a)h)=X\sigma
(a)h=(\varphi _{E}(a)\otimes I)Xh=(\varphi _{E}(a)\otimes I)u(X\otimes h)$.
\end{proof}

The following proposition formalizes in terms of our notation and
perspective an important point made by Skeide in the sentence before Theorem
2.3 of \cite{mS03}. It will help to clarify the relations between our work
here and results in \cite{BS00} and \cite{BBLS04}. To state it, we let $E$
be a $W^{\ast }$-correspondence over the von Neumann algebra $M$ and we let $%
\sigma $ be a faithful normal representation of $M$ on a Hilbert space $H$.
We write $End(E)$ for the collection of all adjointable \emph{bimodule} maps
on $E$. Thus an element of $End(E)$ is an element $T$ in $\mathcal{L}(E)$
that commutes with $\varphi _{E}(M)$; i.e., $End(E)$ is the commutant of $%
\varphi _{E}(M)$ in $\mathcal{L}(E)$. Form the dual correspondence $%
E^{\sigma }$. Then an element $T$ in $End(E)$ \textquotedblleft
induces\textquotedblright\ an element in $End(E^{\sigma })$ that we shall
denote by $T^{\sigma }$. It is defined by the formula $T^{\sigma
}(X)=(T\otimes I)\circ X$, $X\in E^{\sigma }$. Observe that $T^{\sigma }$
makes sense: An $X$ in $E^{\sigma }$ is a map from $H$ to $E\otimes _{\sigma
}H$ and $T\otimes I$ is a map on $E\otimes _{\sigma }H$; so the composition
is a map from $H$ to $E\otimes _{\sigma }H$, i.e, $T^{\sigma }(X)\in
B(H,E\otimes _{\sigma }H)$. To see that $T^{\sigma }(X)\in E^{\sigma }$ for
all $X\in E^{\sigma }$, simply note that since $X\in E^{\sigma }$ and since $%
T\in End(E)$, the desired inclusion results from an easy computation: $%
T^{\sigma }(X)\sigma (a)=(T\otimes I)X\sigma (a)=(T\otimes I)(\varphi
_{E}(a)\otimes I)X=(\varphi _{E}(a)\otimes I)(T\otimes I)X$, for all $a\in M$%
,. Also, it is clear that $T^{\sigma }(Xb)=T^{\sigma }(X)b$ for all $b\in
\sigma (M)^{\prime }$ and that $\left\Vert T^{\sigma }\right\Vert \leq
\left\Vert T\otimes I\right\Vert =\left\Vert T\right\Vert $. Thus, since $%
E^{\sigma }$ is a $W^{\ast }$-correspondence over $\sigma (M)^{\prime }$, we
conclude that $T^{\sigma }$ is automatically adjointable and, therefore,
that $T^{\sigma }\in \mathcal{L}(E^{\sigma })$. The fact that $T^{\sigma }$
lies in $End(E^{\sigma })$ is quite evident since for $b\in \sigma
(M)^{\prime }$ and $X\in E^{\sigma }$, $\varphi _{\sigma }(b)X=(I_{E}\otimes
b)X$ by definition.

\begin{proposition}
\label{dualend}If $E$ is a $W^{\ast }$-correspondence over a von Neumann
algebra $M$ and if $\sigma $ is a faithful normal representation of $M$ on a
Hilbert space $H$, then the map from $End(E)$ to $End(E^{\sigma })$ defined
by the formula $T\rightarrow T^{\sigma }$ is a normal $\ast $-isomorphism
from $End(E)$ onto $End(E^{\sigma })$.
\end{proposition}

\begin{proof}
The only things that really need comment are the facts that the map is
isometric and surjective. Both will be shown with the aid of the Hilbert
space isomorphism $u:E^{\sigma }\otimes _{\iota }H\rightarrow E\otimes
_{\sigma }H$ from Lemma \ref{u}.

By definition of $u$ and $T^{\sigma }$ it is clear that $u(T^{\sigma
}\otimes I_{H})=(T\otimes I_{H})u$. Since the induced representations of $%
\mathcal{L}(E^{\sigma })$ and $\mathcal{L}(E)$ are faithful, since $\iota $
and $\sigma $ are faithful, we conclude that the map $T\rightarrow T^{\sigma
}$ is isometric.

If $S\in End(E^{\sigma })$, let $R=u(S\otimes I_{H})u^{\ast }$. Then $R\in
B(E\otimes _{\sigma }H)$. To show that $R\in \mathcal{L}(E)\otimes I_{H}$,
it suffices to show that $R$ commutes with the produced representation of $%
\sigma (M)^{\prime }$, by Lemma \ref{comm}. But this is clear from the
intertwining properties of $u$: For $b\in \sigma (M)^{\prime }$,
\begin{eqnarray*}
R(I\otimes b) &=&u(S\otimes I)u^{\ast }(I\otimes b)=u((S\otimes I)(\varphi
_{\sigma }(b)\otimes I))u^{\ast } \\
&=&u(S\varphi _{\sigma }(b)\otimes I)u^{\ast }=u(\varphi _{\sigma
}(b)S\otimes I)u^{\ast } \\
&=&u((\varphi _{\sigma }(b)\otimes I)(S\otimes I))u^{\ast }=(I\otimes
b)u(S\otimes I)u^{\ast }=(I\otimes b)R\text{.}
\end{eqnarray*}%
Thus $R=T\otimes I$ for a suitable operator $T\in \mathcal{L}(E)$. To see
that $T\in End(E)$ is another calculation based on Lemma \ref{u}: for $a\in
M $, we have
\begin{eqnarray*}
T\varphi _{E}(a)\otimes I &=&(T\otimes I)(\varphi _{E}(a)\otimes
I)=u(S\otimes I)u^{\ast }(\varphi _{E}(a)\otimes I) \\
&=&u((S\otimes I_{H})(I_{E^{\sigma }}\otimes \sigma (a)))u^{\ast }
=u((S\otimes I_{H})(I_{E^{\sigma }}\otimes \sigma (a)))u^{\ast } \\
&=&u((I_{E^{\sigma }}\otimes \sigma (a))(S\otimes I_{H}))u^{\ast } =(\varphi
_{E}(a)\otimes I)u(S\otimes I)u^{\ast } \\
&=&(\varphi _{E}(a)\otimes I)(T\otimes I)=\varphi _{E}(a)T\otimes I\text{.}
\end{eqnarray*}%
Again, since $\sigma $ is faithful, the induced representation of $\mathcal{L%
}(E)$ is faithful, and we conclude that $T\in End(E)$.
\end{proof}

\begin{lemma}
\label{qT}Let $E$ be a $W^{\ast }$-correspondence over $M$ and let $%
(T,\sigma )$ be a covariant representation of $E$ on a Hilbert space $H$,
where $\sigma $ is faithful. Then

\begin{enumerate}
\item[(i)] There is a projection $q_{T}\in End(E)$ such that
\begin{equation}
q_{T}E\otimes H=[(I_{E}\otimes \sigma (M)^{\prime })\tilde{T}^{\ast }H].
\label{q}
\end{equation}

\item[(ii)] $(I-q_{T})E=\ker T$.

\item[(iii)] $T$ is injective if and only if $q_{T}=I$.
\end{enumerate}
\end{lemma}

\begin{proof}
Let $L$ be the closed subspace of $E\otimes _{\sigma }H$ spanned by the
vectors of the form $(I_{E}\otimes b)\tilde{T}^{\ast }h$ for $b\in \sigma
(M)^{\prime },h\in H$. We also write $[(I_{E}\otimes \sigma (M)^{\prime })%
\tilde{T}^{\ast }H]$ for $L$. Since $L$ is invariant under $I_{E}\otimes
\sigma (M)^{\prime }$, the projection onto it lies in the commutant of $%
I_{E}\otimes \sigma (M)^{\prime }$, which is $\mathcal{L}(E)\otimes I_{H}$.
So we write $q_{T}\otimes I$ for this projection. It follows from the
covariance property of $\tilde{T}$ that $L$ is also invariant under $\varphi
(M)\otimes I_{H}$. Thus $q_{T}\in \mathcal{L}(E)\cap \varphi (M)^{\prime
}=End(E)$, which proves (i). For (ii), note that, given $\xi \in E$, $\xi $
belongs to $(I-q_{T})E$ if and only if for all $h,k\in H$ and all $b\in
\sigma (M)^{\prime }$, $0=\langle \xi \otimes h,(I_{E}\otimes b)\tilde{T}%
^{\ast }k\rangle =\langle \tilde{T}(\xi \otimes bh),k\rangle =\langle T(\xi
)bh,k\rangle $; that is, if and only if $T(\xi )=0$. The third assertion
follows immediately from (ii).
\end{proof}

\begin{definition}
\label{SupportProj}The projection $q_{T}$ associated with a covariant
representation $(T,\sigma )$ of a $W^{\ast }$-correspondence $E$ in Lemma %
\ref{qT} is called the \emph{support projection }of $T$.
\end{definition}

Central to our study is the connection between completely positive maps on a
von Neumann algebra and representations of their Arveson-Stinespring
correspondences. To describe this connection, which was established in \cite%
{MSQMP}, fix a normal, contractive, completely positive map $\Theta $ on a
von Neumann algebra $N$, which we shall assume is represented faithfully on
a Hilbert space $H$. We omit reference to this representation. Form the
Arveson-Stinespring correspondence for $\Theta $, $E_{\Theta }$, which,
recall is a correspondence over $N^{\prime }$. Also, let $W_{\Theta
}:H\rightarrow N\otimes _{\Theta }H$ be the map defined by the equation $%
W_{\Theta }h:=1\otimes h$. Since $\Theta $ is contractive, $W_{\Theta }$ is
a contraction. ($W_{\Theta }$ is an isometry if and only if $\Theta $ is
unital.) A calculation shows that $W_{\Theta }^{\ast }$ is given by the
formula $W_{\Theta }^{\ast }(X\otimes h)=\Theta (X)h$ (see equation 2.1 in
\cite{MSQMP}.) We define $T:E_{\Theta }\rightarrow B(H)$ by the formula
\begin{equation}
T(X):=W_{\Theta }^{\ast }X\text{,}  \label{identityrep}
\end{equation}%
$X\in E_{\Theta }$. If $\sigma $ denotes the identity representation of $%
N^{\prime }$ on $H$, then the pair $(T,\sigma )$ is a representation of $%
E_{\Theta }$ in the sense of Definition \ref{repres}, which is called the
\emph{identity} representation of $E_{\Theta }$. We then have the following
equation, which was proved in Corollary 2.23 of \cite{MSQMP}, showing how to
express $\Theta $ in terms of $(T,\sigma )$:%
\begin{equation}
\Theta (S)=\tilde{T}(I_{E_{\Theta }}\otimes S)\tilde{T}^{\ast }\text{,}
\label{cprep}
\end{equation}%
$S\in N$. Note that $\tilde{T}$ is $W_{\Theta }^{\ast }u$, where $%
u:E_{\Theta }\otimes H\rightarrow M\otimes _{\Theta }H$ is the Hilbert space
isomorphism of Lemma \ref{u} that sends $X\otimes h$ in $E_{\Theta }\otimes
H $ to $Xh$ in $M\otimes _{\Theta }H$. We note, too, that it is possible to
represent all the positive powers of $\Theta $ through similar formulae \cite%
[Theorem 2.24]{MSQMP}.

The following lemma shows that in the representation $(T,\sigma )$ that
arises from a completely positive map $\Theta $, $T$ always is injective.

\begin{lemma}
\label{1dir}Suppose $N$ is a von Neumann algebra acting on the Hilbert space
$H$ and let $\Theta $ be a normal, contractive, completely positive map on $%
N $. Let $E_{\Theta }$ be its Arveson-Stinespring correspondence
(constructed with respect to the identity representation of $N$ on $H$) and
let $(T,\sigma )$ be the identity representation of $E_{\Theta }$. Then

\begin{enumerate}
\item[(i)] $T$ is injective (i.e. its support projection $q_{T}$ is $%
I_{E_{\Theta }}$).

\item[(ii)] $\Theta $ is multiplicative (i.e., an endomorphism) if and only
if $\tilde{T}^{\ast }\tilde{T}=I_{E\otimes H}$.

\item[(iii)] $\Theta $ is unital if and only if $\tilde{T}\tilde{T}^{\ast
}=I $.
\end{enumerate}
\end{lemma}

\begin{proof}
In order to prove that $q_T=I$, fix $X \in E_{\Theta}$ with $X =(I-q_T)X$.
Then, it follows from the definition of $q_T$, that, for every $h,k\in H$
and every $b\in N$, we have $0=\langle \tilde{T}^*k,X \otimes bh \rangle
=\langle k, \tilde{T}(X \otimes bh)\rangle = \langle k, T(X)bh \rangle =
\langle k,W^*X bh \rangle = \langle I\otimes k, X bh \rangle =\langle
I\otimes k, (b\otimes I_H)X h \rangle =\langle b^* \otimes k,X h \rangle$.
But, since $b^*\otimes k$ (for $b\in N$ and $k\in H$) generate $%
N\otimes_{\Theta}H$, we get $X=0$. Thus $q_T=I$.

For (ii), recall that $W_{\Theta }^{\ast }(b\otimes h)=\Theta (b)h$ for $%
b\in N,h\in H$. Then, for $b,c\in N,h,k\in H$, $\langle W_{\Theta }^{\ast
}(b\otimes h),W_{\Theta }^{\ast }(c\otimes k)\rangle =\langle \Theta
(b)h,\Theta (c)k\rangle $. We use the multiplicativity of $\Theta $ to
conclude that this is equal to $\langle h,\Theta (b^{\ast }c)k\rangle
=\langle b\otimes h,c\otimes k\rangle $, showing that $W_{\Theta }$ is an
coisometry. Then, for $X,Y\in E_{\Theta }$, $T(X)^{\ast }T(Y)=X^{\ast
}WW^{\ast }Y=X^{\ast }Y=\sigma (\langle X,Y\rangle )$. Thus $(T,\sigma )$ is
isometric.

Assertion (iii) is immediate from equation (\ref{cprep}).
\end{proof}

\begin{definition}
\label{coisomety}A completely contractive covariant representation $%
(T,\sigma )$ of $E$ is said to be \emph{fully coisometric} if $\tilde{T}$ is
a coisometry; that is, if $\tilde{T}\tilde{T}^{\ast }=I_{H}$.
\end{definition}

The reasons for the terminology regarding a covariant representation $%
(T,\sigma )$, \textquotedblleft isometric\textquotedblright\ and
\textquotedblleft coisometric\textquotedblright , is fairly clear - the
accompanying operator $\tilde{T}$ must be an isometry or a coisometry. The
reason for the adverb \textquotedblleft fully\textquotedblright\ is somewhat
more complicated to explain. \ For this, we refer the reader to Section 5 of
\cite{MSTensor}. It is important for us that the completely positive map $%
\Theta $ is unital if and only if $(T,\sigma )$ is fully coisometric.

Conversely, suppose $E$ is a $W^{\ast }$-correspondence over a von Neumann
algebra $M$ and that $(T,\sigma )$ is a completely contractive covariant
representation of $E$ on a separable Hilbert space $H$, with $\sigma $
faithful. Then for $b\in \sigma (M)^{\prime }$, we define%
\begin{equation}
\Theta _{T}(b)=\tilde{T}(I_{E}\otimes b)\tilde{T}^{\ast }.  \label{thetat}
\end{equation}%
Then by \cite[Proposition 2.21]{MSQMP}, $\Theta _{T}$ is a well defined,
normal, contractive completely positive map on $N:=\sigma (M)^{\prime }$. If
we apply the preceding analysis to $\Theta _{T}$, it is natural to ask how
the correspondence that is produced and the representation of it are related
to $E$ and $(T,\sigma )$. We will show that they are essentially the same,
provided $T$ is injective.

So start with a $W^{\ast }$-correspondence $E$ over a von Neumann algebra $M$
and a completely contractive, covariant representation $(T,\sigma )$ of $E$
on $H$ with faithful $\sigma $. Write $\Theta =\Theta _{T}$ for the normal,
completely positive map defined on $\sigma (M)^{\prime }$ as in equation (%
\ref{thetat}). Applying the discussion above (equations (\ref{etheta}) and (%
\ref{identityrep})) to this $\Theta $, we get a $W^{\ast }$-correspondence $%
E_{\Theta }$ (over $\sigma (M)$) and a covariant representation of this
correspondence, denoted $(T_{\Theta },\sigma _{\Theta })$, where $\sigma
_{\Theta }$ is the identity representation of $\sigma (M)$. Then we have:

\begin{theorem}
\label{Inverse}In the setup just described, there is an isomorphism of
correspondences $w:E_{\Theta }\rightarrow q_{T}E$, where $q_{T}$ is the
support projection of $T$, that carries $(T_{\Theta },\sigma _{\Theta })$ to
the restriction of $(T,\sigma )$ to $q_{T}E$. More precisely, $w$ is a
one-to-one map from $E_{\Theta }$ onto $q_{T}E$ satisfying

\begin{enumerate}
\item[(i)] For $a,b \in M$ and $X\in E_{\Theta}$, $w(\sigma(a)\cdot X
\sigma(b) )=\varphi(a)w(X)b$.

\item[(ii)] For $X,Y\in E_{\Theta }$, $\sigma (\langle w(X),w(Y)\rangle
)=\langle X,Y\rangle $.

\item[(iii)] $\sigma \circ \sigma ^{-1}=\sigma _{\Theta }$ - the identity
representation of $\sigma (M)$ - and $T\circ w=T_{\Theta }$.
\end{enumerate}

In particular, if $T$ is injective, $w$ is an isomorphism of correspondences
from $E_{\Theta}$ onto $E$.
\end{theorem}

\begin{proof}
First, write $N$ for $\sigma (M)^{\prime }$ and let $v:N\otimes _{\Theta
}H\rightarrow E\otimes _{\sigma }H$ be defined by $v(b\otimes h)=(I\otimes b)%
\tilde{T}^{\ast }h$ (and extended by linearity). Since $\langle (I\otimes b)%
\tilde{T}^{\ast }h,(I\otimes c)\tilde{T}^{\ast }k\rangle =\langle h,\Theta
(b^{\ast }c)k\rangle =\langle b\otimes h,c\otimes k\rangle $, for $b,c\in
\sigma (M)^{\prime }$ and $h,k\in H$, the map $v$ is an isometry. It follows
from Lemma~\ref{qT} that the range of $v$ is $q_{T}E\otimes H$. Also note
that, for $b\in \sigma (M)^{\prime }$, $v(b\otimes I)=(I\otimes b)v$ and,
for $a\in M$, $v(I\otimes \sigma (a))=(\varphi (a)\otimes I)v$ (where the
latter equality follows from the covariance property of $\tilde{T}$).

Fix $X\in E_{\Theta }$ and write
\begin{equation*}
\psi (\xi )h=X^{\ast }v^{\ast }(\xi \otimes h)\text{,}
\end{equation*}%
for every $\xi \in E$ and $h\in H$. For $b\in \sigma (M)^{\prime }$, compute
$\psi (\xi )bh=X^{\ast }v^{\ast }(\xi \otimes bh)=X^{\ast }v^{\ast
}(I\otimes b)(\xi \otimes h)=X^{\ast }(b\otimes I)v^{\ast }(\xi \otimes
h)=bX^{\ast }v^{\ast }(\xi \otimes h)$ (where the last equality uses the
covariance property of $X$). Thus $\psi (\xi )\in \sigma (M)$. Also note
that, for $a\in M$, $\psi (\xi a)=\psi (\xi )\sigma (a)$ (as $\xi a\otimes
h=\xi \otimes \sigma (a)h$). Thus the map $\xi \mapsto \sigma ^{-1}(\psi
(\xi ))$ is a (right) module map from $E$ to $M$. It follows from the
selfduality of $E$ that there is a unique element, $w(X)$, in $E$ such that,
for every $\xi \in E$ and $h\in H$,
\begin{equation}
\sigma (\langle w(X),\xi \rangle )h=X^{\ast }v^{\ast }(\xi \otimes h)
\label{w}
\end{equation}%
For $a,b\in M,\xi \in E$ and $h\in H$, we have $\sigma (\langle w(\sigma
(a)\cdot X\sigma (b)),\xi \rangle )h=\sigma (\langle w((I\otimes \sigma
(a))X\sigma (b)),\xi \rangle )h=\sigma (b^{\ast })X^{\ast }(I\otimes \sigma
(a^{\ast }))v^{\ast }(\xi \otimes h)=\sigma (b^{\ast })X^{\ast }v^{\ast
}(\varphi (a^{\ast })\otimes I)(\xi \otimes h)=\sigma (b^{\ast })X^{\ast
}v^{\ast }(\varphi (a^{\ast })\xi \otimes h)=\sigma (b^{\ast })X^{\ast
}v^{\ast }(\varphi (a^{\ast })\xi \otimes h)=\sigma (b^{\ast }\langle
w(X),\varphi (a^{\ast })\xi \rangle )h\newline
=\sigma (\langle \varphi (a)w(X)b,\xi \rangle )h$. Thus $w(\sigma (a)\cdot
X\sigma (b))=\varphi (a)w(X)b$, proving (i).

For $\xi \in E$ we write $L_{\xi}$ for the operator $L_{\xi}: H \rightarrow
E\otimes_{\sigma}H$ mapping $h$ to $\xi \otimes h$. It is easy to check that
$L_{\xi}^*(\zeta \otimes h)=\sigma(\langle\xi,\zeta \rangle)h$ so that
equation (\ref{w}) can be written $L_{w(X)}^*L_{\xi}=X^*v^*L_{\xi}$ for all $%
\xi \in E$. Thus
\begin{equation}  \label{lw}
L_{w(X)}=vX
\end{equation}
for every $X\in E_{\Theta}$.

Fixing $X,Y \in E_{\Theta}$, we have $\sigma(\langle
w(X),w(Y)\rangle)=L_{w(X)}^*L_{w(Y)}=X^*v^*vY=X^*Y$ proving (ii).

Next we show that the image of $w$ is $q_TE$. For this, let $\xi$ in $E$ be
orthogonal to the range of $w$. But then $L_{\xi}^*vX=\sigma(\langle
\xi,w(X) \rangle)=0$ for all $X\in E_{\Theta}$. Since the closed subspace
spanned by the ranges of all $X\in E_{\Theta}$ is all of $N\otimes_{\Theta}H$
(\cite[Lemma 2.10]{MSQMP}) and since the image of $v$ is $q_TE\otimes H$, we
get that $\xi$ is orthogonal to $q_TE$. This shows that the range of $w$
contains $q_TE$ but the argument above can be reversed to show that equality
holds.

It is left to prove part (iii). The first equality is obvious (as $\sigma
_{\Theta }$ is the identity map). For the second, fix $X\in E_{\Theta }$ and
compute, for $h,k\in H$, $\langle W^{\ast }Xh,k\rangle =\langle Xh,I\otimes
k\rangle =\langle vXh,v(I\otimes k)\rangle =\langle w(X)\otimes h,\tilde{T}%
^{\ast }k\rangle =\langle \tilde{T}(w(X)\otimes h),k\rangle =\langle
T(w(X))h,k\rangle .$ Thus $T_{\Theta }(X)=W^{\ast }X=T(w(X))$.
\end{proof}

The following corollary of Theorem \ref{Inverse} is immediate from Lemma \ref%
{1dir} and \cite[Proposition 2.21]{MSQMP}.

\begin{corollary}
Suppose that $M$ is a von Neumann algebra and that $E$ is a $W^{\ast }$%
-correspondence over $M$. If $(T,\sigma )$ is a completely contractive
covariant representation of $E$ on a Hilbert space $H$, where $\sigma $ is
faithful and $T$ is injective, then:

\begin{enumerate}
\item[(i)] $(T,\sigma )$ is isometric if and only if $\Theta _{T}$ is
multiplicative.

\item[(ii)] $(T,\sigma )$ is fully coisometric if and only if $\Theta _{T}$
is unital.
\end{enumerate}
\end{corollary}

\begin{remark}
\label{units}We pause to summarize some of the salient features of our
discussion and to help clarify the relation between between the GNS
correspondence of a $cp$-map and its Arveson-Stinespring correspondence. It
is clear from Example \ref{cpmod} that there is a bijective correspondence
between (contractive) normal $cp$-maps of a von Neumann algebra $M$ and
correspondences $E$ over $M$ endowed with \textquotedblleft
bi-cyclic\textquotedblright\ vectors $\xi $ of norm at most one. By
\textquotedblleft bi-cyclic\textquotedblright\ we mean that the vector
generates $E$ as a bimodule over $M$: $E$ is the closed linear span in the $%
\sigma $-topology of \cite{BDH88} of the set of products $a\xi b$, $a,b\in M$%
. Example \ref{cpmod} shows how a $cp$-map gives rise to such a pair $(E,\xi
)$. And conversely, given a pair $(E,\xi )$, we get a $cp$-map $\Theta $ of $%
M$ via the formula: $\Theta (a)=\langle \xi ,a\xi \rangle $, $a\in M$, and
since $\xi $ is \textquotedblleft bi-cyclic\textquotedblright , $(E,\xi )$
is isomorphic to the GNS correspondence of $\Theta .$ In the literature $%
\Theta $ is sometimes referred to as the coefficient of $\xi $ (cf. \cite%
{AnD90, Mi89, P86}). The coefficient of a vector $\xi $ is a \emph{unital }$%
cp$-map if and only if $\xi $ is a \emph{unit }vector in the sense that $%
\langle \xi ,\xi \rangle =1$, the identity of $M$\footnote{
We learned this very suggestive terminology from \cite[Definition 3.1 et.
seq.]{BS00}}. This is much more than saying that $\xi $ has norm one and, of
course, the notion of a unit vector makes sense in the setting of any
Hilbert module over a unital $C^{\ast }$-algebra. On the other hand, if $E$
is any $W^{\ast }$-correspondence over a von Neumann algebra $M$ and if $%
\sigma $ is a faithful normal representation of $M$ on a Hilbert space $H$,
then the elements of $E^{\sigma }$ of norm at most one are in bijective
correspondence with the contractive covariant representations $(T,\sigma )$
of $M$ through the formula $T\rightarrow \tilde{T}^{\ast }$, as we indicated
in the discussion surrounding equation (\ref{covariance}). As we noted in
Proposition \ref{corres}, $E^{\sigma }$ is a correspondence over $\sigma
(M)^{\prime }$ and the coefficient of the vector $\tilde{T}^{\ast }$
determined by $(T,\sigma )$ is the completely positive map $\Theta _{T}$ on $%
\sigma (M)^{\prime }$ defined in equation (\ref{thetat}). So $\Theta _{T}$
is unital precisely when $\tilde{T}^{\ast }$ is a unit vector in $E^{\sigma }
$, which occurs precisely when $\tilde{T}$, viewed as an operator, is a
coisometry, i.e., precisely when $(T,\sigma )$ is fully coisometric.

Suppose that we are given a covariant representation $(T,\sigma )$ of $E$,
where $\sigma $ is faithful.and suppose that $\tilde{T}^{\ast }\in E^{\sigma
}$ is a \textquotedblleft bi-cyclic\textquotedblright\ vector. Then as we
just observed, $E^{\sigma }$ is the GNS module for the $cp$-map $\Theta _{T}$%
. By Theorem \ref{Inverse}, the Arveson-Stinespring correspondence
associated $\Theta _{T}$ is isomorphic to $q_{T}E$. So, if $\tilde{T}^{\ast
} $ is a bi-cyclic vector in $E^{\sigma }$, then $E$ is (essentially) the
Arveson-Stinespring correspondence for $\Theta _{T}$ if and only if $T$ is
injective.
\end{remark}

\section{Semigroups of Completely Positive Maps, their Product Systems and
their Dilations\label{SemigpsCPMsProdSysDil}}

So far we discussed a single completely positive map and the $W^{\ast }$%
-correspondence associated with it. We now recall and develop further some
basic facts about semigroups of completely positive maps and the product
systems associated with them.

First, we need some terminology.

\begin{definition}
\label{cp0} let $N$ be a von Neumann algebra.

\begin{enumerate}
\item[(i)] A $cp$-semigroup on $N$ is a semigroup $\{\Theta_t : t\geq 0\}$
of contractive, completely positive, normal, linear maps on $N$ such that $%
\Theta_0=id$ and, for $t,s \geq 0$, $\Theta_{t+s}=\Theta_t \circ \Theta_s$.

\item[(ii)] A $cp_0$-semigroup on $N$ is a $cp$-semigroup of unital maps.

\item[(iii)] An $e$-semigroup on $N$ is a $cp$-semigroup of endomorphisms.

\item[(iv)] An $e_0$-semigroup on $N$ is a $cp$-semigroup of unital
endomorphisms.
\end{enumerate}
\end{definition}

Note that we are not making any assumptions, at this stage, about the
continuity properties of these semigroups.

\begin{definition}
\label{productsystem} A family $\{E(t)\}_{t\geq 0}$ of $W^{\ast }$%
-correspondences over a von Neumann algebra $M$ is said to be a product
system if $E(0)=M$ (with the trivial right and left actions of $M$) and if
for every $t,s\geq 0$, there is an isomorphism $U_{t,s}$ (of $W^{\ast }$%
-correspondences) mapping $E(t)\otimes _{M}E(s)$ onto $E(t+s)$ such that $%
U_{t+s,r}(U_{t,s}\otimes I_{E(r)})=U_{t,s+r}(I_{E(t)}\otimes U_{s,r})$ for
every $t,s,r\geq 0$ and such that, for every $t\geq 0$, $U_{t,0}$ and $%
U_{0,t}$ are the right and left actions of $M$ on $E(t)$.

We shall refer to the maps $\{U_{t,s}\}_{t,s\geq 0}$ as the \emph{%
multiplication isomorphisms} of the system. Often they will be suppressed in
calculations.
\end{definition}

Implicit in this definition is the assumption that the $E(t)$ are essential
or unital as left $M$-modules, i.e. that $\varphi _{t}(1)$ is the identity
operator in $\mathcal{L}(E(t))$, for all $t\geq 0$.

\begin{definition}
\label{IsoProdSysts}If $\{E(t)\}_{t\geq 0}$ and $\{F(t)\}_{t\geq 0}$ are two
product systems of $W^{\ast }$-correspondences over a von Neumann algebra $M$%
, then an \emph{isomorphism }from $\{E(t)\}_{t\geq 0}$ to $\{F(t)\}_{t\geq 0}
$ is a family $\gamma =\{\gamma _{t}\}_{t\geq 0}$ of correspondence
isomorphisms, with $\gamma _{t}:E(t)\rightarrow F(t)$, that intertwine the
multiplications on $\{E(t)\}_{t\geq 0}$ and $\{F(t)\}_{t\geq 0}$. That is,
if $\{U_{t,s}^{E}\}_{t,s\geq 0}$ and $\{U_{t,s}^{F}\}_{t,s\geq 0}$ are the
multiplication isomorphisms for $\{E(t)\}_{t\geq 0}$ and $\{F(t)\}_{t\geq 0}$%
, respectively, then $\gamma _{t+s}\circ U_{t,s}^{E}=U_{t,s}^{F}\circ
(\gamma _{t}\otimes \gamma _{s})$, for all $t,s\geq 0$.
\end{definition}

Thus a product system is a bundle of correspondences over the positive
half-axis, $[0,\infty )$, where the fibres can be multiplied and two such
product systems are isomorphic if and only if there is a bundle map between
them consisting of isomorphisms that multiply in a certain sense specified
in Definition \ref{IsoProdSysts}.

In the next section we shall deal with measurable product systems. In the
present discussion we do not assume any Borel structure on the system.

Next we want to define the notion of a covariant representation of a product
system, but first, for the sake of clarity, we state a lemma that shows that
certain operations in our definition make sense and that underpins much of
our analysis in this section. It encapsulates facts used freely in \cite%
{MSQMP}, which in turn derive ultimately from \cite{MSTensor}.

\begin{lemma}
\label{mult0} Let $E$ and $F$ be two $W^{\ast }$-correspondences over a von
Neumann algebra $M$ and let $(T,\sigma )$ and $(S,\sigma )$ be covariant
representations of $E$ and $F$ respectively on a Hilbert space $H$, where $%
\sigma $ is faithful.. Write $\Theta _{T}$ and $\Theta _{S}$ for the
corresponding completely positive maps on $\sigma (M)^{\prime }$ and let $%
\Theta =\Theta _{T}\circ \Theta _{S}$. If we set $R(\xi \otimes \zeta
)=T(\xi )S(\zeta )$ for $\xi \in E,\zeta \in F$, then $R$ may be extended
uniquely to $E\otimes F$ so that the extension, along with $\sigma $ forms a
completely contractive covariant representation $(R,\sigma )$ of $E\otimes
_{M}F$ on $H$. Further, we have $\tilde{R}=\tilde{T}(I_{E}\otimes \tilde{S})$
and $\Theta _{R}=\Theta $.
\end{lemma}

\begin{proof}
One checks easily that $\tilde{R}$ defined by the equation $\tilde{R}:=%
\tilde{T}(I_{E}\otimes \tilde{S})$ is a contractive mapping from $E\otimes
F\otimes H$ to $H$ that satisfies $\tilde{R}\sigma (\cdot )=\varphi
_{E\otimes F}(\cdot )\otimes I_{H}\tilde{R}$, i.e., equation (\ref%
{covariance}) is satisfied. By \cite[Lemmas 3.4-3.6]{MSTensor}, there is a
unique completely contractive bimodule map $R:E\otimes F\rightarrow B(H)$
such that $R(\xi \otimes \eta )h=\tilde{R}(\xi \otimes \eta \otimes h)$. The
rest follows by straightforward calculation.
\end{proof}

\begin{definition}
\label{productrep}Given a product system $\{E(t)\}_{t\geq 0}$ of $W^{\ast }$%
-correspondences over a von Neumann algebra $M$, a \emph{(completely
contractive) covariant representation} of $\{E(t)\}_{t\geq 0}$ on a Hilbert
space $H$ is a family $\{T_{t}\}_{t\geq 0}$ where each $T_{t}$ is a
completely contractive linear map from $E(t)$ to $B(H)$ and $T_{0}$ is a
faithful normal representation of $M$ on $H$ such that each pair $%
(T_{t},T_{0})$ is a completely contractive covariant representation of $%
(E(t),M)$ in the sense of Definition \ref{repres} and such that $%
T_{t}\otimes T_{s}=T_{t+s}$ (after identifying $E(t+s)=E(t)\otimes E(s)$)%
\footnote{%
Note that $T_{t}\otimes T_{s}$ makes sense by the preceding lemma and that
for each $t>0$.}.

Such a covariant representation is called \emph{isometric} (respectively,
\emph{fully coisometric }or \emph{injective}) if each $(T_{t},T_{0})$, $t>0$%
, is isometric (respectively, fully coisometric or injective).
\end{definition}

Suppose that $\{\Theta _{t}\}_{t\geq 0}$ is a $cp$-semigroup on a von
Neumann algebra $N$ and assume that $N$ acts faithfully on a separable
Hilbert space $H$. For every $t\geq 0$ write $\mathcal{E}_{t}$ for the
Arveson-Stinespring correspondence $E_{\Theta _{t}}$ (over $N^{\prime }$)
associated with $\Theta _{t}$ as in equation (\ref{etheta}). That is,
\begin{equation}
\mathcal{E}_{t}=E_{\Theta _{t}}=\{X:H\rightarrow N\otimes _{\Theta
_{t}}H|\;Xb=(b\otimes I)X\;\}.  \label{EThetaT}
\end{equation}%
In general, the family $\{\mathcal{E}_{t}\}_{t\geq 0}$ is \emph{not} a
product system. However, for $t,s\geq 0$, there \emph{is} a
\textquotedblleft coisometric multiplication\textquotedblright\ $m_{t,s}$
mapping $\mathcal{E}_{t}\otimes \mathcal{E}_{s}$ onto $\mathcal{E}_{t+s}$
but it may not be an isometry. (This was observed first by Arveson in \cite%
{Arv97}). The definition of $m_{t,s}$ is a bit involved, but is spelled out
in detail in Section 2 of \cite{MSQMP}, see Proposition 2.12, in particular.
We do not need the details of the definition of $m_{t,s}$, but we do want to
record for future reference the fact that%
\begin{equation}
m_{t,s}^{\ast }\text{ embeds }\mathcal{E}_{t+s}\text{ in }\mathcal{E}%
_{t}\otimes \mathcal{E}_{s}  \label{isoembedding}
\end{equation}%
isometrically as a $W^{\ast }$-correspondence map for all $s,t\geq 0$.

As we showed in \cite{MSQMP}, it is possible to \textquotedblleft
refine\textquotedblright\ the family $\{\mathcal{E}_{t}\}_{t\geq 0}$ in
order to obtain a product system $\{E_{\Theta }(t)\}_{t\geq 0}$ over $%
N^{\prime }$. This process was also carried out in \cite{M03} and a
\textquotedblleft dual\textquotedblright\ process was used in \cite{BS00}.
To describe the process from \cite{MSQMP} briefly, we fix $t>0$ and for any
partition $\mathcal{P}=\{t_{0}=0<t_{1}<\ldots <t_{n-1}<t_{n}=t\}$ of the
interval $[0,t]$ we define $H_{\mathcal{P},t}$ to be $N\otimes _{\Theta
_{t_{1}}}N\otimes _{\Theta _{t_{2}-t_{1}}}N\otimes \cdots \otimes _{\Theta
_{t_{n-1}-t_{n-2}}}N\otimes _{\Theta _{t-t_{n-1}}}H$. (Where the latter
space is defined by successively applying the definition of $N\otimes
_{\Theta }H$). If $\mathcal{P}^{\prime }=\{t_{0}=0<t_{1}<\ldots t_{k}<\tau
<t_{k+1}<\ldots t_{n-1}<t_{n}=t\}$ is a one-point refinement of $\mathcal{P}$
then we can embed $H_{\mathcal{P},t}$ isometrically into $H_{\mathcal{P}%
^{\prime },t}$ by sending $a_{1}\otimes a_{2}\cdots a_{n}\otimes h$ to $%
a_{1}\otimes \cdots a_{k}\otimes a_{k+1}\otimes I\otimes a_{k+2}\cdots
a_{n}\otimes h$. Applying this embedding successively, we can construct an
isometric embedding of $H_{\mathcal{P},t}$ into $H_{\mathcal{P}^{\prime },t}$
whenever $\mathcal{P}^{\prime }$ refines $\mathcal{P}$. We denote this
embedding by $v_{0\mathcal{P}^{\prime },\mathcal{P}}$ and note that $v_{0%
\mathcal{P}^{\prime },\mathcal{P}}$ intertwines the action of $N$. Taking
the direct limit (over the set of all partitions of $[0,t]$ ordered by
refinement) we get the Hilbert space $H_{t}$. Note that in the purely
algebraic situation we are describing, $H_{t}$ is not separable in general.
However, when the construction is carried out using an ultraweakly
continuous semigroup, $\{\Theta _{t}\}_{t\geq 0}$ acting on a countably
decomposable von Neumann algebra, which is represented faithfully on a
separable Hilbert space $H_{0}$, then the space $H_{t}$ will be separable
because it can be seen to be a direct limit over a countable cofinal subset
of the collection of all partitions of $[0,t]$ directed by refinement. The
direct limit is, in fact, a direct limit of left $N$-modules because the $%
v_{0\mathcal{P}^{\prime },\mathcal{P}}$ are $N$-module maps, so we end up
with a (normal) action of $N$ on $H_{t}$. We let $E_{\Theta }(t)$ be the
space of all bounded maps from $H$ to $H_{t}$ that intertwine the actions of
$N$. One can then define on $E_{\Theta }(t)$ a structure of a correspondence
over $N^{\prime }$ (similar to the one defined in Proposition~\ref{corres})
making $E_{\Theta }(t)$ a $W^{\ast }$-correspondence over $N^{\prime }$.
Applying \cite[Lemma 3.2]{MSQMP}, we find that the resulting family $%
\{E_{\Theta }(t)\}_{t\geq 0}$ is a product system. We note that in \cite%
{MSQMP}, it was assumed that the $cp$-semigroup is unital but this
assumption was not used for the construction of $\{E_{\Theta }(t)\}_{t\geq 0}
$.

We also get a (canonical) covariant representation $\{T_{t}\}_{t\geq 0}$ of
this product system (see \cite[Equation (3.1) and Theorem 3.9]{MSQMP}) which
we refer to as the \emph{identity representation} of $\{E_{\Theta
}(t)\}_{t\geq 0}$. It is an analogue of the identity representation of a
single completely positive map defined in equation (\ref{identityrep}).

\begin{remark}
\label{endprod} It is important to note here that if the $cp$-semigroup is a
semigroup of endomorphisms (i.e., an $e$-semigroup), then the family $\{%
\mathcal{E}_{t}\}_{t\geq 0}$ is already a product system so that this
\textquotedblleft refining\textquotedblright\ process is not necessary and,
if performed, yields $E_{\Theta }(t)=\mathcal{E}_{t}$. It is natural to ask
which $cp$-semigroups have this property. One of the consequences of our
analysis is that any $cp$-semigroup that is subordinate to an $e$-semigroup
in a sense that we shall define in Section \ref{SubOrdCPSemi} has this
property, i.e., it already is a product system. See Corollary \ref%
{PersistentProductSyst}.
\end{remark}

Conversely, suppose we are given a product system $\{E(t)\}_{t\geq 0}$ over
a von Neumann algebra $M$ and a covariant representation $\{T_{t}\}_{t\geq
0} $ of it on a Hilbert space $H$. Then we can define a $cp$-semigroup on
the von Neumann algebra $T_{0}(M)^{\prime }$ by
\begin{equation}
\Theta _{t}(b)=\tilde{T}_{t}(I_{E(t)}\otimes b)\tilde{T}_{t}^{\ast
},\;\;b\in T_{0}(M)^{\prime }.  \label{thethat}
\end{equation}%
The fact that this is a semigroup follows easily from the equality $\tilde{T}%
_{t+s}=(I_{E(t)}\otimes \tilde{T}_{s})\tilde{T}_{t}$. For details, see \cite[%
Theorem 3.10]{MSQMP}. We shall assume in the sequel that $T_{0}$ is the
identity representation of $M$ and we shall write $N$ for the commutant of $%
M $. As we explained above, the $cp$-semigroup $\{\Theta _{t}\}_{t\geq 0}$
on $N$ defined via equation (\ref{thethat}) gives rise to a product system $%
\{E_{\Theta }(t)\}_{t\geq 0}$ over $N^{\prime }=M$. Our objective is to show
that $\{E_{\Theta }(t)\}_{t\geq 0}$ is, in a natural way, isomorphic to a
subproduct system of $\{E(t)\}_{t\geq 0}$. That is, we want to prove an
analogue of Theorem~\ref{Inverse}.

For this purpose, we shall use the following lemma, Lemma \ref{mult}, but
first it will be helpful to put some of our discussion into a more general
context.

\begin{remark}
\label{TensorProducts}Let $E$ and $F$ be two $W^{\ast }$-correspondences
over the von Neumann algebra $M$. In general, given $T\in \mathcal{L}(E)$
and $S\in \mathcal{L}(F)$, the formula $\xi \otimes \eta \rightarrow T\xi
\otimes S\eta $ does \emph{not} define an element in $\mathcal{L}(E\otimes
F) $. The formula for the inner product in $E\otimes F$ shows what the
problem is: $\langle T(\xi _{1})\otimes S(\eta _{1}),T(\xi _{1})\otimes
S(\eta _{1})\rangle =\langle S(\eta _{1}),\langle T(\xi _{1}),T(\xi
_{2})\rangle S(\eta _{2})\rangle $. Unless one can \textquotedblleft
slip\textquotedblright\ $\langle T(\xi _{1}),T(\xi _{2})\rangle $ past $S$
and write $\langle T(\xi _{1}),T(\xi _{2})\rangle S(\eta _{2})=S(\langle
T(\xi _{1}),T(\xi _{2})\rangle \eta _{2})$, one will run into difficulties
trying to estimate the norms of certain sums that contain expressions like $%
\langle T(\xi _{1})\otimes S(\eta _{1}),T(\xi _{1})\otimes S(\eta
_{1})\rangle $ in terms of the norms of $S$ and $T$. However, if $S$ lies in
$End(F)$, then one can \textquotedblleft slip\textquotedblright\ $\langle
T(\xi _{1}),T(\xi _{2})\rangle $ by $S$ and it is not hard to see that the
formula $\xi \otimes \eta \rightarrow T\xi \otimes S\eta $ defines an
adjointable operator with norm at most $\left\Vert T\right\Vert \left\Vert
S\right\Vert $. In this case, we denote the element of $\mathcal{L}(E\otimes
F)$ so defined by $T\otimes S$. And of course we have $T\otimes S=(T\otimes
I)(I\otimes S)=(I\otimes S)(T\otimes I)$.

In fact, this observation together with the discussions of induced and
produced representations that we have had all fall under a bigger framework,
still: We can discuss correspondences $E$ from one von Neumann algebra, $M,$
say, to another, say $N$. That is $E$ is a $W^{\ast }$-Hilbert module over $%
N $, endowed with a normal $\ast $-representation $\varphi :M\rightarrow
\mathcal{L}(E)$. If we are given also a $W^{\ast }$-correspondence $F$ from $%
N$ to $P$, then we may form their internal, self-dual tensor product $%
E\otimes F$ to obtain a correspondence from $M$ to $P$. An element $T\in
\mathcal{L}(E)$ and an element $S\in \mathcal{L}(F)$ form an operator,
denoted $T\otimes S$ in $\mathcal{L}(E\otimes F)$ if $S$ lies in the
relative commutant of the left action of $N$ on $F$. In the special case
when $F$ is a Hilbert space, i.e., a $W^{\ast }$-Hilbert module over $%
\mathbb{C}$, then the left action of $N$ is nothing more than a normal
representation, say $\sigma $, of $N$ on $F$. For $T\in M$, $T\otimes I$ is
the induced representation of $M$ determined by $E$ and the left action of $%
M $ on $E$, and, of course, $I\otimes S$, $S\in \sigma (N)^{\prime }$ is the
produced action. We will not need all of this, but the broad context may be
helpful for understanding. A good place to get the rudiments of the theory
as well as connections with other work on correspondences is the paper by
Baillet, Denizeau and Havet \cite{BDH88}. Also, what we are discussing falls
under the rubric of $W^{\ast }$-categories. For these, see \cite{GLR85} and
references therein. And of course, correspondences are ubiquitous in
subfactor theory and in quantum field theory. For an encyclopedic account,
together with an extensive bibliography, see \cite{EK98}.
\end{remark}

\begin{lemma}
\label{mult} With the hypotheses and notation as in Lemma \ref{mult0}, let $%
q_{T}$, $q_{S}$, and $q_{R}$ be the support projections of $T,$ $S$, and $R$
defined in Lemma \ref{qT}. Then $q_{R}\leq q_{T}\otimes q_{S}$.
\end{lemma}

\begin{proof}
Note that $\xi \in E$ lies in $(I-q_{T})E$ if and only if $T(\xi )=0$ (Lemma~%
\ref{qT}). Since $R(\xi \otimes \zeta )=T(\xi )S(\zeta )$, for $\xi \in
E,\zeta \in F$, we get that $I-q_{R}$ dominates both $I_{E}\otimes
(I_{F}-q_{S})=I_{E\otimes F}-(I_{E}\otimes q_{S})$ and $I_{E\otimes
F}-(q_{T}\otimes I_{F})$. It follows that $q_{R}\leq q_{T}\otimes q_{S}$.
\end{proof}

To show that $\{E_{\Theta }(t)\}_{t\geq 0}$ is isomorphic to a subsystem of $%
\{E(t)\}_{t\geq 0}$, we write $\mathcal{E}_{t}$ for $E_{\Theta _{t}}$; as in
equation (\ref{EThetaT}) and we write $q_{t}$ for the support projection $%
q_{T_{t}}$ of $T_{t}$. Theorem~\ref{Inverse} provides an isomorphism of
correspondences
\begin{equation*}
w_{t}:\mathcal{E}_{t}\rightarrow q_{t}E(t)
\end{equation*}%
such that, for $X\in \mathcal{E}_{t}$ and $h\in H$,
\begin{equation}
w_{t}(X)\otimes h=v_{t}\circ X(h)\text{,}  \label{wt}
\end{equation}%
where $v_{t}:N\otimes _{\Theta _{t}}H\rightarrow \mathcal{E}_{t}\otimes H$
is defined by the equation $v_{t}(b\otimes _{\Theta _{t}}h)=(I\otimes b)%
\tilde{T}_{t}^{\ast }h$, $b\otimes _{\Theta _{t}}h\in N\otimes _{\Theta
_{t}}H$. Considering $m_{t,s}^{\ast }$ from equation (\ref{isoembedding}) we
obtain the following diagram in which the vertical arrows are isomorphisms of
$W^{\ast }$-correspondences:%
\begin{equation*}
\begin{array}{ccc}
\mathcal{E}_{t+s} & \overset{m_{t,s}^{\ast }}{\longrightarrow } & \mathcal{E}%
_{t}\otimes \mathcal{E}_{s} \\
\begin{array}{ccc}
w_{t+s} & \downarrow  &
\end{array}
&  &
\begin{array}{ccc}
& \downarrow  & w_{t}\otimes w_{s}%
\end{array}
\\
q_{t+s}E(t+s) &  & q_{t}E(t)\otimes q_{s}E(s)%
\end{array}%
\end{equation*}%
Since each projection $q_{s}$ lies in $End(E(s))$ we may conclude from Lemma %
\ref{mult} that
\begin{equation}
q_{t+s}\leq q_{t}\otimes q_{s},  \label{convexity}
\end{equation}%
for all $s,t\geq 0$, after identifying $E(t+s)$ with $E(t)\otimes E(s)$
using $U_{t,s}$. Thus the isomorphisms $w_{t}$ respect the \textquotedblleft
order\textquotedblright\ on the $\mathcal{E}_{t}$ and the order on the
spaces $q_{t}E(t)$. Or to say things a bit more suggestively, the \emph{%
proto-product system }$\{\mathcal{E}_{t}\}_{t\geq 0}$ is isomorphic to the
image of the \textquotedblleft convex\textquotedblright\ family $%
\{q_{t}\}_{t\geq 0}$ which defines a bundle endomorphism of $\{E(t)\}_{t\geq
0}$.

There is equivalent way to define $E_{\Theta }(t)$, which we use here. Write
$\mathcal{E}(\mathcal{P},t)$ for the space of all bounded maps from $H$ to $%
H_{\mathcal{P},t}$ that intertwine the actions of $N$. Then in a fashion
similar to Proposition~\ref{corres}, $\mathcal{E}(\mathcal{P},t)$ becomes a
correspondence over $N^{\prime }$. Further, whenever $\mathcal{P}^{\prime }$
refines $\mathcal{P}$, the map $v_{0\mathcal{P}^{\prime },\mathcal{P}}$ that
embeds $H_{\mathcal{P},t}$ into $H_{\mathcal{P}^{\prime },t}$ defines, via
composition, a map $v_{\mathcal{P}^{\prime },\mathcal{P}}:\mathcal{E}(%
\mathcal{P},t)\rightarrow \mathcal{E}(\mathcal{P}^{\prime },t)$, which is an
$N^{\prime }$-correspondence embedding. By \cite[Lemma 3.1]{MSQMP} we
conclude that
\begin{equation}
E_{\Theta }(t)\simeq \lim_{\rightarrow }(\mathcal{E}(\mathcal{P},t),v_{%
\mathcal{P}^{\prime },\mathcal{P}}).  \label{dirlim}
\end{equation}%
The reason this isomorphism is important for our analysis is that we may
apply \cite[Proposition 2.12]{MSQMP} to construct an isomorphism of
correspondences
\begin{equation}
\Psi (\mathcal{P},t):\mathcal{E}(\mathcal{P},t)\rightarrow \mathcal{E}%
_{t-t_{n-1}}\otimes \cdots \otimes \mathcal{E}_{t_{1}},  \label{PsiP}
\end{equation}%
for each partition $\mathcal{P}=\{0=t_{0}<t_{1}\cdots <t_{n}=t\}$. Thanks to
the discussion in the preceding paragraph, we may identify each space $%
\mathcal{E}_{t_{k}-t_{k-1}}$ with the range of the projection $%
q_{t_{k}-t_{k-1}}$ in $\mathcal{L}(E(t_{k}-t_{k-1}))$ and, using the maps $%
U_{t,s}$ from Definition \ref{productsystem}, we may build a projection in $%
\mathcal{L}(E(t))$, denoted $q(\mathcal{P},t)=q_{t-t_{n-1}}\otimes
q_{t_{n-1}-t_{n-2}}\otimes \cdots \otimes q_{t_{1}}$, and an isomorphism of $%
W^{\ast }$-correspondences $\Xi (\mathcal{P},t)$ from $\mathcal{E}%
_{t-t_{n-1}}\otimes \cdots \otimes \mathcal{E}_{t_{1}}$ onto $q(\mathcal{P}%
,t)E(t)$ so that if $\mathcal{P}^{\prime }$ refines $\mathcal{P}$ then $q(%
\mathcal{P},t)\leq q(\mathcal{P}^{\prime },t)$ in $\mathcal{L}(E(t))$ and so
that the following diagram commutes%
\begin{equation*}
\begin{array}{ccc}
\mathcal{E}(\mathcal{P},t) & \overset{v_{\mathcal{P}^{\prime },\mathcal{P}}}{%
\longrightarrow } & \mathcal{E}(\mathcal{P}^{\prime },t) \\
\downarrow  &  & \downarrow  \\
q(\mathcal{P},t)E(t) & \overset{\iota _{\mathcal{P}^{\prime },\mathcal{P}}}{%
\longrightarrow } & q(\mathcal{P}^{\prime },t)E(t)%
\end{array}%
\end{equation*}%
where the left vertical arrow is $\Xi (\mathcal{P},t)\Psi (\mathcal{P},t),$
the right vertical arrow is $\Xi (\mathcal{P}^{\prime },t)\Psi (\mathcal{P}%
^{\prime },t)$ and where $\iota _{\mathcal{P}^{\prime },\mathcal{P}}$ is the
inclusion map. The fact that $q(\mathcal{P},t)\leq q(\mathcal{P}^{\prime },t)
$ in $\mathcal{L}(E(t))$ may be seen easily by viewing $\mathcal{P}^{\prime }
$ as a succession of one-point refinements of $\mathcal{P}$ and applying the
`convexity' relation (\ref{convexity}), which is a consequence of Lemma \ref%
{mult}: If $\mathcal{P}=\{t_{0}=0<t_{1}<\ldots <t_{n-1}<t_{n}=t\}$ and if $%
\mathcal{P}^{\prime }=\{t_{0}=0<t_{1}<\ldots t_{k}<\tau <t_{k+1}<\ldots
t_{n-1}<t_{n}=t\}$, then $q(\mathcal{P},t)=q_{t-t_{n-1}}\cdots \otimes
q_{t_{k+1}-t_{k}}\otimes \cdots \otimes q_{t_{1}}\leq q_{t-t_{n-1}}\cdots
\otimes q_{t_{k+1}-\tau }\otimes q_{\tau -t_{k}}\otimes \cdots \otimes
q_{t_{1}}=q(\mathcal{P}^{\prime },t)$.) Consequently, if we define $%
\overline{q}_{t}$ by the formula%
\begin{equation}
\overline{q}_{t}:=\sup_{\mathcal{P}}q(\mathcal{P},t)  \label{supqt}
\end{equation}%
where the supremum is taken over all partitions $\mathcal{P}$ of the
interval $[0,t]$, then $\overline{q}_{t}$ is a projection in $\mathcal{L}%
(E(t))$ such that $E_{\Theta }(t)=\overline{q}_{t}E(t)$. Further, with the
aid of the isomorphisms $U_{t,s}$ that define the product system $%
\{E(t)\}_{t\geq 0}$, we conclude that $\overline{q}_{t+s}=\overline{q}%
_{t}\otimes \overline{q}_{s}$.

We summarize our discussion of the last three paragraphs as

\begin{theorem}
\label{subsys}Let $\{E(t)\}_{t\geq 0}$ be a product system over a von
Neumann algebra $M$ and let $\{T_{t}\}_{t\geq 0}$ be a covariant
representation of $\{E(t)\}_{t\geq 0}$ on a Hilbert space $H$ such that $%
T_{0}$ is injective. Let $\Theta =\{\Theta _{t}\}_{t\geq 0}$ be the $cp$%
-semigroup defined by $\{T_{t}\}_{t\geq 0}$ through equation (\ref{thethat})
and let $\{E_{\Theta }(t)\}_{t\geq 0}$ be the product system generated by $%
\Theta $, as was just described. If $q_{t}$ is the support projection of $%
T_{t}$, for each $t$, and if $\overline{q}_{t}$ is defined by equation (\ref%
{supqt}) then $\overline{q}_{t}$ is a projection in $End(E(t))$ with $%
q_{t}\leq \overline{q}_{t}$ for all $t$, $\overline{q}_{t+s}=\overline{q}%
_{t}\otimes \overline{q}_{s}$, for all $s,t\geq 0$, and $\{E_{\Theta
}(t)\}_{t\geq 0}$ is isomorphic to the product system $\{\overline{q}%
_{t}E(t)\}_{t\geq 0}$ - a subproduct system of $\{E(t)\}_{t\geq 0}$.
Further, if $q_{t+s}=q_{t}\otimes q_{s}$ for all $t,s\geq 0$, then $q_{t}=%
\overline{q}_{t}$ for all $t\geq 0$.
\end{theorem}

Note that the last assertion is clear since if $q_{t+s}=q_{t}\otimes q_{s}$
for all $t,s\geq 0$, then $q(\mathcal{P},t)=q_{t}$ for all partitions $%
\mathcal{P}$ of $[0,t]$. So by definition of $\overline{q}_{t}$ as $\sup_{%
\mathcal{P}}q(\mathcal{P},t)$ in equation (\ref{supqt}), the assertion
follows.

\section{Measurable product systems\label{MeasProdSysts}}

In this section we develop the theory of Borel product systems of $W^{\ast }$%
-correspondences and relate them to quantum Markov semigroups. Our
formulation was inspired by Effros's analysis of Borel fields of von Neumann
algebras in \cite{E65}, which seems ready-made to allow us to perform the
analysis we want to undertake. We would like to note that alternate
approaches to our discussion in this section may be possible. For example,
Borel and continuity structures have long played a central role in the
structure of Fell bundles, which are product systems over groups (see \cite%
{FD88}). For the purpose of studying generalizations of the tensor operator
algebras from \cite{MSTensor}, the first author generalized Fell bundles
over groups to Fell bundles over groupoids in \cite{pM01} and showed how to
define a topological product system over a partial order in a topological
groupoid. In \cite{iH04}, Hirshberg gave a definition of a Borel structure
on a system of $C^{\ast }$-correspondences that is close to ours and in \cite%
{mS03a} Skeide attached a continuity structure to an $E_{0}$-semigroup
defined on $\mathcal{L}(E)$, where $E$ is a Hilbert $C^{\ast }$-module over
an auxiliary $C^{\ast }$-algebra. One can get Borel structures from
topologies of course and $W^{\ast }$-correspondences from $C^{\ast }$%
-correspondences, but we shall not delve into the exact relations here,
since they do not contribute to our immediate purposes.

The following definition is modeled on the definition of a Borel field of
Hilbert spaces. (See the \textquotedblleft Borel\textquotedblright\ version
of Definition 1 and Proposition 3 in \cite[Part II, Chapter 1]{Dix}). In it
and throughout the remainder of this paper, when we talk about the Borel
structure on a von Neumann algebra $M$, say, then we mean the Borel
structure generated by the weak operator topology. \ This is the same as the
Borel structure generated by the strong-, ultraweak- or the ultrastrong
operator topology and so, in particular, it is independent of any faithful
representation of $M$.

We want to emphasize that from now on all of our von Neumann algebras are
countably decomposable and all Hilbert spaces under consideration are
separable. These assumptions are essential for most of our analysis, in
particular for dealing with Borel product systems, which are, \emph{inter
alia}, Borel families of $W^{\ast }$-modules in the sense of the following
definition.

\begin{definition}
\label{BorelWmod}Let $M$ be a (countably decomposable) von Neumann algebra.
For every $t\in \lbrack 0,\infty )$, let $E(t)$ be a $W^{\ast }$-module over
$M$. We shall say the the family $\{E(t)\}_{t\geq 0}$ is a \emph{Borel family%
} (or is \emph{measurable}) if there is a countable family of cross sections
$\{f_{n}\}_{n\geq 0}$ (that is, each $f_{n}$ is defined on $[0,\infty )$ and
$f_{n}(t)\in E(t)$ for $t\in \lbrack 0,\infty )$) satisfying:

\begin{enumerate}
\item[(i)] For every $n,m$, the $M$-valued function $t\mapsto \langle
f_{n}(t),f_{m}(t)\rangle $ is Borel measurable.

\item[(ii)] For every $t\in \lbrack 0,\infty )$, the $W^{\ast }$-submodule
of $E(t)$ generated by $\{f_{n}(t)\}_{n\geq 0}$ is $E(t)$.
\end{enumerate}

Such a family of cross sections $\{f_{n}\}_{n\geq 0}$ will be called a \emph{%
defining family }of cross sections for the family $\{E(t)\}_{t\geq 0}$.
\end{definition}

Given a Borel family of $W^{\ast }$-modules $\{E(t)\}_{t\geq 0}$ and a
defining family of cross sections $\{f_{n}\}_{n\geq 0}$, we set $[0,\infty
)\star E:=\{(t,\xi ):\;\xi \in E(t)\}$ and let $\pi :[0,\infty )\star
E\rightarrow \lbrack 0,\infty )$ be the projection onto the first factor, $%
\pi (t,\xi )=t$. Then algebraically, $[0,\infty )\star E$ together with $\pi
$ forms a bundle over $[0,\infty )$. Further, the sections $\{f_{n}\}_{n\geq
0}$ define maps $\{f_{n}^{\prime }\}_{n\geq 0}$ of this bundle to $M$ via
the equation $f_{n}^{\prime }(t,\xi )=\langle f_{n}(t),\xi \rangle $. The
\emph{Borel structure on }$[0,\infty )\star E$ is defined to be the smallest
Borel structure making the map $\pi $ and each of the maps $f_{n}^{\prime }$
Borel.

It is immediate, therefore, that a section $t\mapsto (t,g(t))\in \lbrack
0,\infty )\star E$ is a Borel map if and only if its composition with each $%
f_{n}^{\prime }$ is Borel; that is, if and only if the map to $M$, $t\mapsto
\langle f_{n}(t),g(t)\rangle $ is Borel for each $n$.

\begin{definition}
\label{Borelproduct}A family $\{E(t)\}_{t\geq 0}$ of $W^{\ast }$%
-correspondences over a (countably decomposable) von Neumann algebra $M$ is
said to be a \emph{Borel family of }$W^{\ast }$\emph{-correspondences }(or a
\emph{measurable family of }$W^{\ast }$-\emph{correspondences}) if

\begin{enumerate}
\item[(i)] It is measurable as a family of $W^{\ast }$-modules in the sense
of Definition \ref{BorelWmod}.

\item[(ii)] The map from $M\times ([0,\infty )\star E)$ to $[0,\infty )\star
E$ sending $(a,(t,\xi ))$ to $(t,\varphi _{t}(a)\xi )$ is a Borel map.
\end{enumerate}

The family $\{E(t)\}_{t\geq 0}$ is called a \emph{measurable product system}
(\emph{of }$W^{\ast }$\emph{-correspondences}) if, in addition, it satisfies:

\begin{enumerate}
\item[(iii)] For every $t,s \geq 0$ there is an isomorphism $U_{t,s}$ (of $%
W^*$-correspondences) mapping $E(t) \otimes_M E(s)$ onto $E(t+s)$ such that $%
U_{t+s,r}(U_{t,s}\otimes I_{E(r)})=U_{t,s+r}(I_{E(t)}\otimes U_{s,r})$ for
every $t,s,r \geq 0$.

\item[(iv)] The family $\{U_{s,t}\}_{s,t\geq 0}$ is Borel in the sense that
if $t\mapsto \xi (t)$ is a Borel section and if $\eta $ lies in $E(s)$, then
the sections $l(\eta )\xi $ and $r(\eta )\xi $, defined by the equations
\begin{equation*}
(l(\eta )\xi )(t)=U_{s,t-s}(\eta \otimes \xi (t-s))
\end{equation*}%
\begin{equation*}
(r(\eta )\xi )(t)=U_{t-s,s}(\xi (t-s)\otimes \eta )
\end{equation*}%
are Borel (where $\xi (t-s)$ is understood to be $0$ if $t<s$).
\end{enumerate}
\end{definition}

\begin{remark}

\begin{enumerate}
\item[(1)] Often we shall suppress the maps $\varphi _{t}$ and $U_{t,s}$
when working with a measurable product system of $W^{\ast }$-correspondences.

\item[(2)] It follows from part (ii) of the definition that, whenever $%
t\mapsto a_{t}$ is a Borel map from $[0,\infty )$ to $M$ and $t\mapsto \xi
(t)$ is a Borel section of $\{E(t)\}_{t\geq 0}$, the map $t\mapsto a_{t}\xi
(t)$ $(=\varphi _{t}(a_{t})\xi (t))$ is also a Borel section.

\item[(3)] We do not know if the condition (iv) in the definition of a
measurable product system is equivalent to the measurability of the map from
$([0,\infty )\star E)\times ([0,\infty )\star E)$ to $[0,\infty )\star E$
induced by the family $\{U_{s,t}\}_{s,t\geq 0}$.
\end{enumerate}
\end{remark}

Suppose that $\{E(t)\}_{t\geq 0}$ is a Borel family of $W^{\ast }$-modules
over $M$ and that $\{f_{n}\}_{n\geq 0}$ is a defining family of cross
sections. We also fix a (faithful) representation $\sigma $ of $M$ on a
Hilbert space $H$ and set $H(t)=E(t)\otimes _{\sigma }H$. Then the family of
sections $\{f_{n}(t)\otimes e_{m}\}_{n,m\geq 0}$, where $\{e_{m}\}_{m\geq 0}$
is a fixed orthonormal basis of $H$, makes $\{H(t)\}_{t\geq 0}$ a Borel
family of Hilbert spaces. It then follows from \cite[Proposition 1 in
section II.1.4]{Dix} that there is a fixed Hilbert space $H_{0}$ and unitary
operators $v_{t}$ from $H_{0}$ onto $H(t)$ such that, for every $h\in H$,
the section $t\mapsto v_{t}h$ is Borel. (Note that we can always assume that
each $H(t)$ is infinite dimensional by choosing $H$ to be infinite
dimensional at the outset.) Thus, the measurable Hilbert bundle determined
by the $H(t)$ is trivialized by the unitary family $\{v_{t}\}_{t\geq 0}$. It
follows that given a Borel section $t\mapsto \xi (t)\in E(t)$ and vectors $%
k\in H$ and $h\in H_{0}$, the map $t\mapsto \langle L_{\xi (t)}^{\ast
}v_{t}h,k\rangle =\langle v_{t}h,\xi (t)\otimes k\rangle $ is Borel.

\begin{lemma}
Let $\{E(t)\}_{t\geq 0}$ be a Borel family of $W^{\ast }$-modules and let $%
g_{1}$ and $g_{2}$ be two Borel sections. Then the map $t\mapsto \langle
g_{1}(t),g_{2}(t)\rangle $ is a Borel map.
\end{lemma}

\begin{proof}
For a Borel family of Hilbert spaces, this result is known (\cite[%
Proposition 4, Chapter 1, Part II]{Dix}) . It is possible to prove the
present lemma using similar methods; that is, by applying a Gram-Schmidt
process in the setting of $W^{\ast }$-modules. But it is easier simply to
take vectors $x_{1},x_{2}\in H$ (where $M$ is represented faithfully on $H$%
), note that $\langle g_{1}(t)\otimes x_{1},g_{2}(t)\otimes x_{2}\rangle
=\langle x_{1},\langle g_{1}(t),g_{2}(t)\rangle x_{2}\rangle $, and then use
the result about sections of a Borel family of Hilbert spaces.
\end{proof}

Continuing with our discussion, define an isomorphism $\Phi _{t}:\mathcal{L}%
(E(t))\rightarrow B(H_{0})$ by setting $\Phi _{t}(T)=v_{t}^{\ast }(T\otimes
I_{H})v_{t}$ (for $t\geq 0$) where $\{v_{t}\}_{t\geq 0}$ is a trivializing
family of unitaries as above. Since $\{f_{n}\}_{n\geq 0}$ is a defining
family of cross sections of $\{E(t)\}_{t\geq 0}$, we see that if $h,k$ are
vectors in $H_{0}$, then the map $t\mapsto \langle \Phi _{t}(f_{n}(t)\otimes
f_{m}(t)^{\ast })h,k\rangle =\langle L_{f_{n}(t)}^{\ast
}v_{t}h,L_{f_{m}(t)}^{\ast }v_{t}k\rangle $ is a Borel map for every $n$ and
$m$. Note also, that for a fixed $t\geq 0$, $\Phi _{t}(\mathcal{L}(E(t)))$
is the von Neumann algebra generated by the operators $\Phi
_{t}(f_{n}(t)\otimes f_{m}(t)^{\ast })$. Thus, $\{\Phi _{t}(\mathcal{L}%
(E(t)))\}_{t\geq 0}$ is a measurable family of von Neumann algebras in $%
B(H_{0})$ in the sense of \cite[Definition II.3.2.1]{Dix}.

Now assume that in addition to being a Borel family of $W^{\ast }$-modules, $%
\{E(t)\}_{t\geq 0}$ is a Borel family of $W^{\ast }$-correspondences over $M$%
. Then for $a\in M$ and $n,m\in \mathbb{N}$, the map $t\mapsto
af_{n}(t)\otimes e_{m}$ is Borel and so is the map $t\mapsto \langle
(\varphi _{t}(a)\otimes I)v_{t}h_{0},f_{n}(t)\otimes e_{m}\rangle =\langle
v_{t}h_{0},a^{\ast }f_{n}(t)\otimes e_{m}\rangle $ (for $h_{0}\in H_{0}$).
It follows that, for $a\in M$, the map that sends $t$ to the matrix $\left(
\begin{array}{cc}
\Phi _{t}(\varphi _{t}(a)) & 0 \\
0 & a%
\end{array}%
\right) $ is a Borel map. If, for each $t$, $B_{t}$ is the von Neumann
algebra generated by these matrices as $a$ runs over $M$, we get a Borel
family of von Neumann subalgebras of $B(H_{0}\oplus H)$.

By a result of Effros \cite[Theorem 3]{E65} the field $t\mapsto
B_{t}^{\prime }$ is also a Borel field of von Neumann algebras. Thus, there
is a family $\{R_{n}(t)\}_{n\geq 0}$ of Borel maps into $B(H_{0}\oplus H)$
such that, for every $t$, the set $\{R_{n}(t)\}_{n\geq 0}$ generates the von
Neumann algebra $B_{t}^{\prime }$. Letting $P$ and $P_{0}$ be the
projections of $H_{0}\oplus H$ onto $H$ and $H_{0}$, respectively, we find
that operators in the right upper corner (i.e operators in $%
P_{0}B_{t}^{\prime }P$) are precisely the operators $T:H\rightarrow H_{0}$
that satisfy the equation $v_{t}^{\ast }(\varphi _{t}(a)\otimes I)v_{t}T=Ta$
for all $a\in M$. It follows from the definition of the $\sigma $-dual, $%
E(t)^{\sigma }$, (where $\sigma $ is the fixed representation of $M$ on $H$)
that these are precisely the operators $T$ such that $v_{t}T$ lies in $%
E(t)^{\sigma }$.

Hence the family $\{E(t)^{\sigma }\}_{t\geq 0}$ is a Borel family of $%
W^{\ast }$-modules over $M^{\prime }$, where $\{v_{t}P_{0}R_{n}(t)P\}_{n\geq
0}$ is a defining family of cross sections for the Borel structure.

\begin{theorem}
\label{dualprodsys}Let $M$ be a countably decomposable von Neumann algebra,
let $\sigma $ be a faithful representation of $M$ on a Hilbert space $H$ and
let $\{E(t)\}_{t\geq 0}$ be a measurable product system over $M$. Then the
system of duals $\{E(t)^{\sigma }\}_{t\geq 0}$ is a measurable product
system over $\sigma (M)^{\prime }$.
\end{theorem}

\begin{proof}
We have already shown that $\{E(t)^{\sigma }\}_{t\geq 0}$ is a Borel family
of $W^{\ast }$-modules. So we need to verify conditions (ii)-(iv) in
Definition~\ref{Borelproduct}. To verify condition (ii), fix a Borel map $%
t\mapsto b_{t}\in M^{\prime }$ and a Borel section $t\mapsto \eta _{t}$ of $%
\{E(t)^{\sigma }\}_{t\geq 0}$. Note first that for every $h,k\in H$ and
every pair of Borel sections $t\mapsto \xi _{t}$ and $t\mapsto \zeta _{t}$
of $\{E(t)\}_{t\geq 0}$, we have $\langle (I_{E(t)}\otimes b_{t}^{\ast })\xi
_{t}\otimes h,\zeta _{t}\otimes k\rangle =\langle b_{t}^{\ast }h,\langle \xi
_{t},\zeta _{t}\rangle k\rangle $. Thus $t\mapsto I_{E(t)}\otimes
b_{t}^{\ast }$ is a Borel section of the family $\{B(H_{t})\}_{t\geq 0}$
and, consequently, $t\mapsto v_{t}^{\ast }(I_{E(t)}\otimes b_{t}^{\ast
})v_{t}$ is a Borel map into $B(H_{0})$. Since $t\mapsto \eta _{t}$ is a
Borel section of $\{E(t)^{\sigma }\}_{t\geq 0}$, $t\mapsto v_{t}^{\ast }\eta
_{t}$ is a Borel map into $B(H,H_{0})$ and it follows that, for every $h\in
H $ and $h_{0}\in H_{0}$, the map $t\mapsto \langle v_{t}^{\ast
}(I_{E(t)}\otimes b_{t})\eta _{t}h,h_{0}\rangle =\langle v_{t}^{\ast }\eta
_{t}h,v_{t}^{\ast }(I_{E(t)}\otimes b_{t}^{\ast })v_{t}h_{0}\rangle $ is
Borel. Thus $t\mapsto v_{t}^{\ast }(I_{E(t)}\otimes b_{t})\eta _{t}$ is a
Borel map into $B(H,H_{0})$. Consequently, $t\mapsto (I_{E(t)}\otimes
b_{t})\eta _{t}=b_{t}\cdot \eta _{t}$ is a Borel section of $\{E(t)^{\sigma
}\}_{t\geq 0}$, which verifies (ii).

For (iii), we note that if we write $U_{s,t}$ for the system of isomorphisms
associated with the product system $\{E(t)\}_{t\geq 0}$, the system of
isomorphisms for the dual system, denoted $\{W_{s,t}\}_{s,t\geq 0}$, are
defined by the equation
\begin{equation}
W_{s,t}(\zeta \otimes \eta )=(U_{t,s}\otimes I_{H})(I_{E(t)}\otimes \zeta
)\eta   \label{dualisos}
\end{equation}%
where $\zeta \in E(s)^{\sigma }$ and $\eta \in E(t)^{\sigma }$. The fact
that they satisfy (iii) follows from \cite[Lemma 3.7]{MSHardy}.

To verify (iv), fix the following objects: $s\geq 0,\zeta \in E(s)^{\sigma }$%
, $k$ in $H$, $k_{0}$ in $H_{0}$ and a Borel section $\{\eta _{t}\}_{t\geq
0} $ of $\{E(t)^{\sigma }\}_{t\geq 0}$. Also, let $\{f_{n}(t)\}_{n\geq 0}$
be a defining family of Borel sections for Borel structure on $%
\{E(t)\}_{t\geq 0}$. From the fact that $\{U_{s,t}\}_{s,t\geq 0}$ is
measurable (in the sense of Definition~\ref{Borelproduct} (iv)) we see that
for every $\xi \in E(s)$ and $h\in H$, the map $t\mapsto \langle
U_{t-s,s}(f_{n}(t-s)\otimes \xi )\otimes h,v_{t}k_{0}\rangle =\langle
(U_{t-s,s}\otimes I_{H})(f_{n}(t-s)\otimes (\xi \otimes
h)),v_{t}k_{0}\rangle $ is Borel. Let $\{e_{m}\}_{m\geq 0}$ be a fixed
orthonormal basis for $H$. Then $\zeta (e_{m})$ lies in $E(s)\otimes H$ and,
thus the map $t\mapsto \langle (U_{t-s,s}\otimes I_{H})(f_{n}(t-s)\otimes
\zeta (e_{m})),v_{t}k_{0}\rangle $ is Borel (for all $n,m$). The latter
expression can also be written as $\langle f_{n}(t-s)\otimes
e_{m},(I_{E(t-s)}\otimes \zeta ^{\ast })(U_{t-s,s}^{\ast }\otimes
I_{H})v_{t}k_{0}\rangle $ and so the map $t\mapsto (I_{E(t-s)}\otimes \zeta
^{\ast })(U_{t-s,s}^{\ast }\otimes I_{H})v_{t}k_{0}$ is a Borel section of $%
\{E(t-s)\otimes H\}_{t\geq 0}$ . Since $t\mapsto \eta _{t-s}k$ is also a
Borel section, their inner product yields a Borel map. Thus $t\mapsto
\langle (U_{t-s,s}\otimes I_{H})(I_{E(t-s)}\otimes \zeta )\eta
_{t-s}k,v_{t}k_{0}\rangle =\langle v_{t}^{\ast }W_{s,t-s}(\zeta \otimes \eta
_{t-s})k,k_{0}\rangle $ is a Borel map. This proves the measurability of $%
l(\zeta )$. The proof of the measurability of $r(\zeta )$ is similar and we
omit it.
\end{proof}

It will be convenient to adopt the following terminology and notation which
are borrowed from \cite{Arv03}.

\begin{definition}
\label{CP0}Let $N$ be a von Neumann algebra.

\begin{enumerate}
\item[(i)] A $CP$-semigroup on $N$ is a semigroup $\{\Theta_t : t\geq 0\}$
of contractive, completely positive, normal, linear maps on $N$ such that,
for every $a\in M$, the map $t \mapsto \Theta_t(a)$ is ultraweakly
continuous.

\item[(ii)] A $CP_0$-semigroup on $N$ is a $CP$-semigroup of unital maps.

\item[(iii)] An $E$-semigroup on $N$ is a $CP$-semigroup of endomorphisms.

\item[(iv)] An $E_{0}$-semigroup on $N$ is a $CP$-semigroup of unital
endomorphisms.
\end{enumerate}
\end{definition}

Of course, a $CP$-semigroup is what we called a quantum Markov semigroup at
the outset of this paper.

The following example shows how to get a Borel product system of
correspondences over a von Neumann algebra $M$ from an $E_{0}$-semigroup on $%
M$. In the setting of discrete product systems of $C^{\ast }$%
-correspondences, the algebraic aspects of this product system are due to
Fowler \cite{nF02}. In all but name, it occurs in the context of semigroups
of endomorphisms of $C^{\ast }$-algebras in Khoshkam and Skandalis's paper
\cite{KS97}. Skeide, too, discusses it without reference to the Borel
structure in \cite{mS03}.

\begin{example}
\label{endom}Let $M$ be a countably decomposable von Neumann algebra acting
on a Hilbert space $H$ and write $\sigma $ for the inclusion of $M$ into $%
B(H)$. Let $\{\alpha _{t}\}_{t\geq 0}$ be an $E$-semigroup acting on $M$ and
let $_{\alpha _{t}}M$ be the correspondence associated with $\alpha _{t}$ as
in Example~\ref{cpmod}. The family $t\mapsto {}_{\alpha _{t}}M$ is a product
system with multiplication isomorphisms $U_{t,s}(a\otimes b)=\alpha _{s}(a)b$
for $a\in {}_{\alpha _{t}}M$ and $b\in {}_{\alpha _{s}}M$, and it is easy to
verify that it is measurable. If we let $E_{\alpha }(t)=E_{\alpha
_{t}}^{\sigma }$ be the $\sigma $-dual correspondence of $_{\alpha _{t}}M$
as in Example \ref{cpdual}, it then follows from Theorem \ref{dualprodsys}
that $t\mapsto E_{\alpha }(t)$ is also a measurable product system, with
multiplication isomorphisms $W_{s,t}(S\otimes T)=ST$.
\end{example}

The first part of the next theorem can be deduced from Theorem~\ref{subsys}.
However, since the proof here is simpler and needs no \textquotedblleft
direct limit arguments'', we prefer to spell it out.

\begin{theorem}
\label{isomprod}Let $\{E(t)\}_{t\geq 0}$ be a product system over $M$ and
let $\{V_{t}\}_{t\geq 0}$ be an isometric covariant representation of it.
For every $t>0$ write $\alpha _{t}$ for the normal endomorphism of $%
V_{0}(M)^{\prime }$ defined as in equation (\ref{thethat}); that is, $\alpha
_{t}(b)=\tilde{V}_{t}(I_{E(t)}\otimes b)\tilde{V}_{t}^{\ast }$ for $b\in
V_{0}(M)^{\prime }$. Then there is an isomorphism $\gamma =\{\gamma
_{t}\}_{t\geq 0}$ of product systems from $\{E_{\alpha }(t)\}_{t\geq 0}$
onto $\{E(t)\}_{t\geq 0}$.

Moreover, if $\{\alpha _{t}\}_{t\geq 0}$ is an $E$-semigroup (that is, the
map $t\mapsto \alpha _{t}(b)$ is ultraweakly continuous for every $b\in
V_{0}(M)^{\prime }$), then $\{E(t)\}_{t\geq 0}$ is a measurable product
system.
\end{theorem}

\begin{proof}
We write $N$ for $V_{0}(M)^{\prime }$. As usual, we shall suppress reference
to the identity representation of $N$. We shall use the isomorphisms
constructed in Theorem~\ref{Inverse}, i.e. the $W^{\ast }$-correspondence
isomorphisms $w_{t}$ from $E_{\alpha }(t)=E_{\alpha _{t}}^{\sigma }$ onto $%
E(t)$. There we proved the equality $L_{w_{t}(X_{t})}=v_{t}X_{t}$ (see
equation (\ref{lw})), where $v_{t}:N\otimes _{\alpha _{t}}H\rightarrow
E(t)\otimes H$ maps $b\otimes h$ to $(I\otimes b)\tilde{V}_{t}^{\ast }h$.
Since the maps $\alpha _{t}$ are endomorphisms, we may identify $N\otimes
_{\alpha _{t}}H$ with $H$ through the unitary map $u_{t}:N\otimes _{\alpha
_{t}}H\rightarrow H$ sending $b\otimes h$ to $\alpha _{t}(b)h$. As is
indicated at the end of Example~\ref{cpdual}, $E_{\alpha }(t)$ is then
identified with the space of all bounded operators in $B(H)$ that intertwine
the identity representation of $N$ and $\alpha _{t}$.

Given such an operator $Y_{t}$ in $B(H)$, we get from equation (\ref{lw})
the equation $L_{w_{t}(u_{t}^{\ast }Y_{t})}=v_{t}u_{t}^{\ast }Y_{t}$. But $%
u_{t}v_{t}^{\ast }=\tilde{V}_{t}$ (since $u_{t}v_{t}^{\ast }(I\otimes b)%
\tilde{V}_{t}^{\ast }h=u_{t}(b\otimes h)=\alpha _{t}(b)h=\tilde{V}%
_{t}(I\otimes b)\tilde{V}_{t}^{\ast }h$) and, thus, $L_{w_{t}(u_{t}^{\ast
}Y_{t})}=\tilde{V}_{t}^{\ast }Y_{t}$. Write $\gamma
_{t}(Y_{t}):=w_{t}(u_{t}^{\ast }Y_{t})$. Then each $\gamma _{t}$ is an
isomorphism of $E_{\alpha }(t)$ onto $E(t)$ and we conclude that $L_{\gamma
_{t}(Y_{t})}=\tilde{V}_{t}^{\ast }Y_{t}$. After applying $\tilde{V}_{t}$, we
find that
\begin{equation}
V_{t}(\gamma _{t}(Y_{t}))=\alpha _{t}(I)Y_{t},\;\;Y_{t}\in E_{\alpha }(t).
\label{vg}
\end{equation}%
For $\xi \in E(t),a\in N$ and $h\in H$, $\alpha _{t}(a)V_{t}(\xi )h=\tilde{V}%
_{t}(I_{E(t)}\otimes a)\tilde{V}_{t}^{\ast }\tilde{V}_{t}(\xi \otimes h)=%
\tilde{V}_{t}(\xi \otimes ah)=V_{t}(\xi )ah$. Thus $V_{t}(\xi )$ lies in $%
E_{\alpha }(t)$. Applying equation (\ref{vg}) with $V_{t}(\xi )$ in place of
$Y_{t}$ (and noting that $\alpha _{t}(I)V_{t}(\xi )=V_{t}(\xi )$), we get $%
V_{t}(\gamma _{t}(V_{t}(\xi )))=V_{t}(\xi )$. Since $V_{t}$ is injective, we
have 
\begin{equation}
\gamma _{t}(V_{t}(\xi ))=\xi ,\;\;\xi \in E(t).  \label{gv}
\end{equation}%
In order to show that $\{\gamma _{t}\}_{t\geq 0}$ gives an isomorphism of
product systems, we need to check that, for $Y_{t}\in E_{\alpha }(t)$ and $%
Y_{s}\in E_{\alpha }(s)$, we have $\gamma _{t+s}(U_{t,s}^{\alpha
}(Y_{t}\otimes Y_{s}))=U_{t,s}(\gamma _{t}(Y_{t})\otimes \gamma _{s}(Y_{s}))$
(where $U$ and $U^{\alpha }$ denote the \textquotedblleft
multiplication\textquotedblright\ maps on $\{E(t)\}_{t\geq 0}$ and $%
\{E_{\alpha }(t)\}_{t\geq 0}$ respectively). Recall from Example~\ref{endom}
that $U_{t,s}^{\alpha }(Y_{t}\otimes Y_{s})=Y_{t}Y_{s}$. Using equation (\ref%
{vg}) and the fact that $V_{t}(\gamma _{t}(Y_{t}))$ lies in $E_{\alpha }(t)$%
, we compute:
\begin{multline*}
V_{t+s}(\gamma _{t+s}(U_{t,s}^{\alpha }(Y_{t}\otimes Y_{s})))=V_{t+s}(\gamma
_{t+s}(Y_{t}Y_{s}))=\alpha _{t+s}(I)Y_{t}Y_{s}=\alpha _{t}(\alpha
_{s}(I))\alpha _{t}(I)Y_{t}Y_{s} \\
=\alpha _{t}(\alpha _{s}(I))V_{t}(\gamma _{t}(Y_{t}))Y_{s}=V_{t}(\gamma
_{t}(Y_{t}))V_{s}(\gamma _{s}(Y_{s})) \\
=V_{t+s}(U_{t,s}(\gamma _{t}(Y_{t})\otimes \gamma _{s}(Y_{s}))).
\end{multline*}%
Applying $\gamma _{t+s}$ to this equality (and using the equation $\gamma
_{t+s}(V_{t+s}(\xi ))=\xi $, $\xi \in E(t+s)$) we get the required result.

If $\{\alpha _{t}\}_{t\geq 0}$ is an $E$-semigroup, we already know that $%
\{E_{\alpha }(t)\}_{t\geq 0}$ is a measurable product system (Example~\ref%
{endom}) and we can use $\gamma =\{\gamma _{t}\}_{t\geq 0}$ to
\textquotedblleft carry'' the Borel structure to $\{E(t)\}_{t\geq 0}$ ; that
is, if $\{Y^{(n)}\}$ is a countable family of sections defining the Borel
structure in $\{E_{\alpha }(t)\}_{t\geq 0}$ , then $\{\gamma (Y^{(n)})\}$
define a Borel structure on $\{E(t)\}_{t\geq 0}$.
\end{proof}

Our next goal is to show that the product system for a $CP$-semigroup \ is
measurable. For this purpose, we will use the preceding Theorem, Theorem \ref%
{isomprod}, and a dilation result. The dilation result, which we state here,
was proved for a fully coisometric covariant representation in \cite[Theorem
3.7]{MSQMP}. The general case was recently proved in \cite[Theorem 1.1]{mS06}%
.

\begin{theorem}
\label{dilrep} (\cite{MSQMP}, \cite{mS06}) Let $\{E(t)\}_{t\geq 0}$ be a
product system over a von Neumann algebra $M$ and let $\{T_{t}\}_{t\geq 0}$
be a fully coisometric covariant representation of the product system on a
Hilbert space $H$. Then there is another Hilbert space $K$, an isometry $%
u_{0}$ mapping $H$ into $K$, and fully coisometric, isometric covariant
representation $\{V_{t}\}_{t\geq 0}$ of $\{E(t)\}_{t\geq 0}$ on $K$ such that

\begin{enumerate}
\item[(1)] $u_0^*V_t(\xi)u_0=T_t(\xi)$ for all $\xi\in E(t), t\geq 0$.

\item[(2)] For all $\xi \in E(t), t\geq 0$, $V_t(\xi)^*$ leaves $u_0(H)$
invariant.

\item[(3)] The smallest subspace of $K$ containing $u_0(H)$ and reducing
each $V_t(\xi)$ is all of $K$.
\end{enumerate}
\end{theorem}

The following theorem was proved, for the case $M=B(H)$ and the product
system involved is a product system of Hilbert spaces, by Markiewicz in \cite%
[Theorem 3.9]{M03}, using different methods.

\begin{theorem}
\label{meas} Let $\{\Theta _{t}\}_{t\geq 0}$ be a $CP$-semigroup on a von
Neumann algebra $M$ and let $\{E_{\Theta }(t)\}_{t\geq 0}$ be the product
system associated to this semigroup as in \cite{MSQMP} (see the discussion
in the previous section). Then $\{E_{\Theta }(t)\}_{t\geq 0}$ is a
measurable product system.
\end{theorem}

\begin{proof}
As we briefly described in Section \ref{SemigpsCPMsProdSysDil}, one can
associate with the $CP$-semigroup a product system $\{E_{\Theta
}(t)\}_{t\geq 0}$ and a representation of this product system on a Hilbert
space. (This is the identity representation mentioned in Section \ref%
{SemigpsCPMsProdSysDil}). We can use Theorem~\ref{dilrep} to dilate this
representation to an isometric representation $\{V_{t}\}_{t\geq 0}$ of $%
\{E_{\Theta }(t)\}_{t\geq 0}$ on $K$. Now we apply Theorem~\ref{isomprod}.
\end{proof}

\begin{definition}
Let $\{E(t)\}_{t\geq 0}$ be a measurable product system over $M$ and let $%
\{T_{t}\}_{t\geq 0}$ be a covariant representation of $\{E(t)\}_{t\geq 0}$
on a Hilbert space $H$. Then the representation is said to be \emph{%
measurable} if, whenever $t\mapsto \xi _{t}$ is a Borel section of the
product system and $h,k\in H$, the map $t\mapsto \langle T_{t}(\xi
_{t})h,k\rangle $ is Borel.
\end{definition}

As we see in the following theorem, such representations (satisfying an
additional condition) give rise to $CP$-semigroups on $T_0(M)^{\prime}$

\begin{theorem}
\label{cont} Let $\{E(t)\}_{t\geq 0}$ be a measurable product system over $M$
and let $\{T_{t}\}_{t\geq 0}$ be a measurable covariant representation of
this product system on $H$. Write $\Theta _{t}$ for the completely positive
map on $T_{0}(M)^{\prime }$ defined by
\begin{equation}
\Theta _{t}(b)=\tilde{T}_{t}(I_{E(t)}\otimes b)\tilde{T}_{t}^{\ast }.
\label{contthetat}
\end{equation}%
Then $\{\Theta _{t}\}_{t\geq 0}$ is a semigroup of normal, contractive,
completely positive maps on $T_{0}(M)^{\prime }$ and, for every $b\in
T_{0}(M)^{\prime }$, the map $t\mapsto \Theta _{t}(b)$ is ultraweakly
continuous on $(0,\infty )$.

If, in addition,
\begin{equation}  \label{0cont}
\cap\{ker(\Theta_t) | t>0 \}=\{0\},
\end{equation}
then the map is also continuous at $t=0$.
\end{theorem}

\begin{proof}
Everything but the continuity can be found in \cite[Theorem 3.10]{MSQMP}. So
we attend to that. Fix $b\in T_{0}(M)^{\prime }$. For every measurable
section $t\mapsto \xi _{t}$ of $\{E(t)\}_{t\geq 0}$ and every $h,k\in H$,
the map that sends $t$ to $\langle (I_{E(t)}\otimes b)\tilde{T}_{t}^{\ast
}h,\xi _{t}\otimes k\rangle =\langle h,T_{t}(\xi _{t})bk\rangle $ is
measurable, since $\{T_{t}\}_{t\geq 0}$ is a measurable representation. Thus
$t\mapsto (I_{E(t)}\otimes b)\tilde{T}_{t}^{\ast }h$ is a Borel section of $%
\{E(t)\otimes H\}_{t\geq 0}$ and so is $t\mapsto \tilde{T}_{t}^{\ast }k$.
Forming the inner product, we conclude that the map $t\mapsto \langle \Theta
_{t}(b)h,k\rangle =\langle (I_{E(t)}\otimes b)\tilde{T}_{t}^{\ast }h,\tilde{T%
}_{t}^{\ast }k\rangle $ is measurable. Write $\mathcal{R}$ for $%
T_{0}(M)^{\prime }$. Then it follows that the function $t\mapsto \langle
\omega ,\Theta _{t}(b)\rangle $ is measurable for every $b\in \mathcal{R}$
and every $\omega \in \mathcal{R}_{\ast }$ (where $\langle \cdot ,\cdot
\rangle $ is the duality pairing of $\mathcal{R}$ and $\mathcal{R}_{\ast }$
). Since each $\Theta _{t}$ is normal, we can write $\Psi _{t}$ for the
pre-adjoint of $\Theta _{t}$, i.e. $\Psi _{t}(\omega )=\omega \circ \Theta
_{t}$ for all $\omega \in \mathcal{R}_{\ast }$.

Since $\mathcal{R}_*$ is separable, \cite[Theorem 3.5.3]{HP57} implies that $%
t \mapsto \Psi_t(\omega)$ is strongly measurable as an $\mathcal{R}_*$%
-valued function. Thus, in the terminology of \cite[ Chap. 10]{HP57}, $%
\{\Psi_t\}_{t\geq 0}$ is a strongly measurable semigroup of linear maps on $%
\mathcal{R}_*$. But then, \cite[Theorem 10.2.3]{HP57} shows that $t \mapsto
\Psi_t$ is strongly continuous on $(0,\infty)$; i.e., for each $\omega \in
\mathcal{R}_*$, the $\mathcal{R}_*$-valued function on $(0,\infty)$, $%
t\mapsto \Psi_t(\omega)=\omega \circ \Theta_t$, is continuous with respect
to the norm topology on $\mathcal{R}_*$. This proves the required continuity
on $(0,\infty)$.

To prove the continuity at $t=0$, assuming condition (\ref{0cont}), we write
$\tilde{\mathcal{R}}_{\ast }$ for the closed linear span $\vee \{\Psi _{t}(%
\mathcal{R}_{\ast })|t>0\}$ and note that, if $\tilde{\mathcal{R}}_{\ast
}\neq \mathcal{R}_{\ast }$, there is a non zero operator $R$ in $\mathcal{R}$
such that $\langle \omega ,R\rangle =0$ for all $\omega \in \tilde{\mathcal{R%
}}_{\ast }$. This means that, for all $t>0$ and all $\omega \in \mathcal{R}%
_{\ast }$, $\langle \omega ,\Theta _{t}(R)\rangle =\langle \Psi _{t}(\omega
),R\rangle =0$. Thus $R$ lies in the kernels of all $\Theta _{t}$, $t>0$,
contradicting condition (\ref{0cont}). It follows that $\tilde{\mathcal{R}}%
_{\ast }=\mathcal{R}_{\ast }$. We can now appeal to \cite[Theorem 10.5.5]%
{HP57} to conclude that $\lim_{t\rightarrow 0+}||\Psi _{t}(\omega )-\omega
||=0$. Consequently, for all $\omega \in \mathcal{R}_{\ast }$ and $R\in
\mathcal{R}$ we see that $\langle \omega ,\Theta _{t}(R)\rangle =\langle
\Psi _{t}(\omega ),R\rangle \rightarrow \langle \omega ,R\rangle $ as $%
t\rightarrow 0+$, which is what we wanted to prove.
\end{proof}

\begin{remark}
\label{weakstrong}Note that the arguments of the proof of the theorem show
that, if $\{\Theta _{t}\}_{t\geq 0}$ is a $CP$-semigroup and $(\Theta
_{t})_{\ast }$ is the pre-adjoint of $\Theta _{t}$, then the map $t\mapsto
(\Theta _{t})_{\ast }(\omega )=\omega \circ \Theta _{t}$ is norm continuous
for all $\omega $ in the predual.
\end{remark}

\section{Subordinate CP-semigroups\label{SubOrdCPSemi}}

As an application of our analysis of measurable product systems and their
relation to $CP$-semigroups, we want to study the notion of
\textquotedblleft subordination\textquotedblright\ for these semigroups. Our
analysis rests on a result of Arveson \cite[Theorem 1.4.2]{Arv69}. In order
to state it, we require some notation. Given a $C^{\ast }$-algebra $B$ and a
Hilbert space $H$, we write $CP(B,H)$ for the set of all completely positive
linear maps from $B$ into $B(H)$. There is a natural partial ordering on
this set defined by $\psi \leq \phi $ if $\phi -\psi $ is completely
positive and, for $\phi \in CP(B,H)$, we write $[0,\phi ]:=\{\psi \in
CP(B,H):\psi \leq \phi \}$. Given $\phi \in CP(B,H)$, we write $\phi
(b)=V^{\ast }\pi (b)V$ for the unique (up to unitary equivalence) minimal
Stinespring dilation of $\phi $, where $\pi $ is the Stinespring
representation of $B$ on the space Hilbert space $K$, $V$ is the contractive
Stinespring \textquotedblleft embedding\textquotedblright\ of $H$ into $K$,
and where $K$ is spanned by $\pi (B)V(H)$. For every $T\in B(K)$ that
commutes with the image of $\pi $ we write $\phi _{T}$ for the map $\phi
_{T}(b)=V^{\ast }T\pi (b)V$. Then we have the following

\begin{proposition}
\label{subarv} (\cite{Arv69}) The map $T\mapsto \phi_T$ is an affine order
isomorphism of the partially ordered convex set of operators $\{T\in
\pi(B)^{\prime}: 0\leq T \leq I \}$ onto $[0,\phi]$.
\end{proposition}

\begin{definition}
\label{subordinate} Let $\Theta$ and $\Psi$ be two normal, contractive,
completely positive maps on the von Neumann algebra $N$. We say that $\Psi$
is subordinate to $\Theta$ (and write $\Psi\leq \Theta$) if $\Theta - \Psi$
is completely positive.

Similarly, if $\{\Theta _{t}\}_{t\geq 0}$ and $\{\Psi _{t}\}_{t\geq 0}$ are $%
cp$-semigroups on $N$, then we say that $\{\Psi _{t}\}_{t\geq 0}$ is
subordinate to $\{\Theta _{t}\}_{t\geq 0}$ if $\Psi _{t}\leq \Theta _{t}$
for all $t\geq 0$.
\end{definition}

We shall adapt Proposition \ref{subarv} to deal with normal completely
positive maps on a given von Neumann algebra $N$, mapping $N$ into itself.
Our completely positive maps are associated with covariant representations
of $W^{\ast }$-correspondences. So let $E$ be such a correspondence (over
the von Neumann algebra $M$, say) and let $(T,\sigma )$ be a covariant
representation of $E$ with faithful $\sigma $. The resulting completely
positive map is $\Theta (b)=\tilde{T}(I_{E}\otimes b)\tilde{T}^{\ast }$ for $%
b\in \sigma (M)^{\prime }$.

Write $N$ for $\sigma (M)^{\prime }$, recall that the Stinespring
representation space $K$ is $N\otimes _{\sigma \circ \Theta }H$ and that, in
the proof of Theorem~\ref{Inverse}, we considered the isometry $v$ mapping $%
K=N\otimes _{\sigma \circ \Theta }H$ onto the subspace $L=[(I\otimes N)%
\tilde{T}^{\ast }H]$ of $E\otimes _{\sigma }H$ (by mapping $b\otimes h$ to $%
(b\otimes I)\tilde{T}^{\ast }h$). Recall also, from Lemma~\ref{qT}, that
there is a projection $q_{T}\in End(E)$ such that $L$ is the range of $%
q_{T}\otimes I_{H}$.

Using $v$ to carry the Stinespring representation $\pi$ from $K$ to $L$, we
get a representation $\tau$ on $L$ defined by $\tau(x)(I \otimes y)\tilde{T}%
^*h=(I\otimes xy)\tilde{T}^*h$ for $x,y \in \sigma(M)^{\prime},\;h\in H$.
Since the commutant of $I \otimes \sigma(M)^{\prime}$ (in $B(E \otimes H)$)
is $\mathcal{L}(E)\otimes I$ (Lemma~\ref{comm}), the commutant of the image
of $\tau$ (in $B(L)$) is the algebra $q_T\mathcal{L}(E)q_T\otimes I_H$. It
follows from Proposition~\ref{subarv} that there is an order isomorphism
between the completely positive maps $\phi: \sigma(M)^{\prime}\rightarrow
B(H)$ with $\phi \leq \Theta$ and positive contractive operators in $q_T%
\mathcal{L}(E)q_T\otimes I_H$. In fact, if $c\in q_T\mathcal{L}(E)q_T$ then
the associated completely positive map is $\Theta_c(x)=\tilde{T}(c\otimes x)%
\tilde{T}^*$.

It follows from this expression that every such map is automatically a
normal map. In general, the image of $\Theta _{c}$ is not contained in $%
\sigma (M)^{\prime }$. To see when it is contained there, fix $a\in M,x\in
\sigma (M)^{\prime }$ and $c$ as above and compute.
\begin{equation}
\Theta _{c}(x)\sigma (a)=\tilde{T}(c\otimes x)\tilde{T}^{\ast }\sigma (a)=%
\tilde{T}(c\otimes x)(\varphi (a)\otimes I)\tilde{T}^{\ast }=\tilde{T}%
(c\varphi (a)\otimes x)\tilde{T}^{\ast }
\end{equation}%
and, similarly,
\begin{equation}
\sigma (a)\Theta _{c}(x)=\tilde{T}(\varphi (a)c\otimes x)\tilde{T}^{\ast }.
\end{equation}%
It follows that the image of $\Theta _{c}$ is contained in $\sigma
(M)^{\prime }$ if and only if $c\in q_{T}(End(E))q_{T}$. Summarizing this
discussion, we have the following.

\begin{proposition}
\label{sub}Let $\Theta $ be the normal contractive completely positive map
associated to a covariant representation $(T,\sigma )$ (with faithful $%
\sigma $) of the $W^{\ast }$-correspondence $E$ over $M$. Then there is an
order isomorphism $c\mapsto \Theta _{c}$ from the set of all operators $c\in
q_{T}(End(E))q_{T}$ that satisfy $0\leq c\leq I$ onto the set of all normal
contractive completely positive maps on $\sigma (M)^{\prime }$ that are
subordinate to $\Theta $.

Given $c$ the map $\Theta _{c}$ is written
\begin{equation*}
\Theta _{c}(x)=\tilde{T}(c\otimes x)\tilde{T}^{\ast },\;\;x\in \sigma
(M)^{\prime }.
\end{equation*}
\end{proposition}

\begin{corollary}
\label{subord}

\begin{enumerate}
\item[(1)] Let $E$ be a $W^{\ast }$-correspondence over a von Neumann
algebra $M$ and let $(T,\sigma )$ and $(S,\rho )$ be two covariant
representations of $E$ with faithful $\sigma $ and $\rho $. Write $q_{T}$
and $q_{S}$ for the support projections of these representations in $End(E)$
defined through equation (\ref{q}). Also, write $\Theta _{T}$ and $\Theta
_{S}$ for the completely positive maps associated with the representations.

If $q_{T}=q_{S}$ then there is an order isomorphism between the normal,
completely positive, maps of $\sigma (M)^{\prime }$ that are subordinate to $%
\Theta _{T}$ and the normal completely positive maps of $\rho (M)^{\prime }$
that are subordinate to $\Theta _{S}$.

\item[(2)] Suppose $\Theta $ is a normal, unital, completely positive map of
a von Neumann algebra $N$ and let $\alpha $ be its minimal endomorphic
dilation acting on the von Neumann algebra $R$. Then there is an order
isomorphism between the normal completely positive maps on $N$ subordinate
to $\Theta $ and those that are subordinate to $\alpha $ on $R$.
\end{enumerate}
\end{corollary}

\begin{proof}
Part (1) follows immediately from Proposition~\ref{sub}. For part (2),
recall the construction of $\alpha$ from \cite{MSQMP}. Given such $\Theta$,
we fix a Hilbert space $H$ on which $N$ acts and we define the
correspondence $E_{\Theta}$ as in (\ref{etheta}). $E_{\Theta}$ is a
correspondence over $N^{\prime}$. We write $\sigma$ for the identity
representation of $N^{\prime}$ on $H$ and define $T:E_{\Theta} \rightarrow
B(H)$ by $T(\eta)=W^*\eta$ where $W:H \rightarrow N\otimes_{\Theta}H$ is the
operator mapping $h\in H$ to $I\otimes h$. The pair $(T,\sigma)$ is then a
fully coisometric covariant representation of $E_{\Theta}$ on $H$ (referred
to as the identity representation). Applying \cite[Corollary 5.21]{MSTensor}
(or \cite[Theorem 2.18]{MSQMP}), we can dilate $(T,\sigma)$ to a covariant
representation $(V,\rho)$ on a Hilbert space $K$ that is both isometric and
fully coisometric. Writing $\alpha(S)=\tilde{V}(I_{E_{\Theta}}\otimes S)%
\tilde{V}^*$ we get a unital endomorphism of $\rho(N^{\prime})^{\prime}$
that dilates $\Theta$. This is the minimal endomorphic dilation of $\Theta$.

The proof of (2) will be complete (using part (1)) once we show that $%
q_{T}=q_{V}$. Since $(V,\rho )$ is isometric it is clear that $q_{V}=I$. The
fact that $q_{T}=I$ was proved in Lemma~\ref{1dir}. Thus $q_{V}=q_{T}=I$.
\end{proof}

The referee has kindly noted that Proposition \ref{sub} and Corollary \ref%
{subord} may be derived from Paschke's Proposition 5.4 in \cite{Pa73} using
duality.

\begin{definition}
\label{reducedcocycles}Let $\{E(t)\}_{t\geq 0}$ a measurable product system
over a von Neumann algebra $M$, let $\{U_{s,t}\}_{s,t\geq 0}$ be the family
of multiplication maps that identify $E(s+t)$ with $E(s)\otimes E(t)$, $%
s,t\geq 0$ and let $q=\{q_{t}\}$ be a family projections such that $q_{t}$
lies in $End(E(t))$ for each $t$, and such that for every $s,t$, $%
q_{s+t}\leq U_{t,s}(q_{t}\otimes q_{s})U_{t,s}^{\ast }$. A family $%
c=\{c_{t}\}_{t\geq 0}$ of sections of $\{\mathcal{L}(E(t))\}_{t\geq 0}$ will
be called a \emph{reduced cocycle for the product system }$\{E(t)\}_{t\geq
0} $ (relative to $q$) if $c$ satisfies

\begin{enumerate}
\item[(Cq1)] For each $t\geq 0$, $c_{t}\in q_{t}(End(E(t)))q_{t}$ and $0\leq
c_{t}\leq I$.

\item[(Cq2)] $c_{0}=I$ and, for $t,s\geq 0$, $c_{t+s}=q_{t+s}U_{s,t}(c_{s}%
\otimes c_{t})U_{t,s}^{\ast }q_{t+s}$.

\item[(Cq3)] $c$ is a Borel section of the Borel family $\{\mathcal{L}%
(E(t))\}_{t\geq 0}$.

\item[(Cq4)] For every $s>0$ and every $\xi ,\zeta \in q_{s}E(s)$, $%
\lim_{t\rightarrow 0}\langle U_{(s-t),t}(I_{E(s-t)}\otimes c_{t})\xi ,\zeta
\rangle =\langle \xi ,\zeta \rangle $ in the ultraweak topology of $M$.
\end{enumerate}

For a given family $q$, as above, the collection of all reduced cocycles for
the product system $\{E(t)\}_{t\geq 0}$ will be denoted $\mathcal{C}_{q}(E)$.

If $q_{t}=I_{E(t)}$ for every $t$, then we call $c$ simply a \emph{cocycle}
for $\{E(t)\}_{t\geq 0}$ and we will write $\mathcal{C}(E)$ for the
collection of all cocycles for $\{E(t)\}_{t\geq 0}$.
\end{definition}

Very often we will omit the multiplication maps and simply write $%
c_{t+s}=c_{s}\otimes c_{t}$. We note, too, that the notion of a (reduced)
cocycle for a product system is different from the more familiar notion of a
cocycle for a semigroup. \ However, as we shall see in Proposition \ref%
{cocycleequal}, there is a relation between the two notions.

As we have seen in our construction of a product system from a $cp$%
-semigroup, families of projections $q=\{q_{t}\}$ satisfying the hypothesis $%
q_{t+s}\leq U_{t,s}(q_{t}\otimes q_{s})U_{t,s}^{\ast }$ for all $s,t\geq 0$
in Definition \ref{reducedcocycles} arise as the support projections of a
representation. This will play an important role in what follows. To support
our computations, we begin with a lemma.

\begin{lemma}
\label{multbis} Let $E$ and $F$ be two $W^{\ast }$-correspondences over a
von Neumann algebra $M$ and let $(T,\sigma )$ and $(S,\sigma )$ be covariant
representations of $E$ and $F$, respectively, on $H$, where $\sigma $ is
assumed to be faithful. Let $(R,\sigma )$ be the representation of $E\otimes
F$ defined by the formula $R(\xi \otimes \zeta )=T(\xi )S(\zeta )$ for $\xi
\in E,\zeta \in F$ and described in Lemma \ref{mult0}. Write $\Theta _{T}$, $%
\Theta _{S}$ and $\Theta _{R}$ for the corresponding completely positive
maps on $\sigma (M)^{\prime }$. Suppose $\Psi $ and $\Phi $ are normal,
completely positive maps on $\sigma (M)^{\prime }$ such that $\Psi \leq
\Theta _{T}$ and $\Phi \leq \Theta _{S}$, and let $c_{1}\in End(E)$ and $%
c_{2}\in End(F)$ be the positive elements such that $\Psi =(\Theta
_{T})_{c_{1}}$ and $\Phi =(\Theta _{T})_{c_{2}}$. Then $\Psi \circ \Phi
=\Theta _{c}$ where $c=q_{R}(c_{1}\otimes c_{2})q_{R}$, where $q_{R}$ is the
support projection of $R$.
\end{lemma}

\begin{proof}
The proof is a computation: $\Psi \circ \Phi (b)=\Psi (\tilde{S}%
(c_{1}\otimes b)\tilde{S}^{\ast })=\tilde{T}(c_{2}\otimes (\tilde{S}%
(c_{1}\otimes b)\tilde{S}^{\ast }))\tilde{T}^{\ast }=\tilde{R}(c_{1}\otimes
c_{2}\otimes b)\tilde{R}^{\ast }=\tilde{R}(q_{R}(c_{1}\otimes
c_{2})q_{R}\otimes b)\tilde{R}^{\ast }$, for $b\in \sigma (M)^{\prime }$.
\end{proof}

\begin{theorem}
\label{subordinatetheorem}Let $\{E(t)\}_{t\geq 0}$ be a measurable product
system and let $\{T_{t}\}_{t\geq 0}$ be a measurable covariant
representation of it. Write $q=\{q_{t}\}$ where $q_{t}=q_{T_{t}}$ is the
support projection of $T_{t}$ for every $t$, and suppose that the semigroup $%
\Theta =\{\Theta _{t}\}_{t\geq 0}$ of completely positive maps of $%
T_{0}(M)^{\prime }$ induced by $\{T_{t}\}_{t\geq 0}$ (defined via equation (%
\ref{thethat})) is a $CP$-semigroup. Then there is an order isomorphism $%
c\mapsto \Theta _{c}$ from the set $\mathcal{C}_{q}(E)$ of reduced cocycles
for $\{E(t)\}_{t\geq 0}$, relative to $q$, onto the set of all $CP$%
-semigroups $\Psi =\{\Psi _{t}\}_{t\geq 0}$ of $T_{0}(M)^{\prime }$ that are
subordinate to $\Theta $.

The semigroup associated to $c\in \mathcal{C}_{q}(E)$ is given by the
equation%
\begin{equation}
(\Theta _{c})_{t}(b)=\tilde{T}_{t}(c_{t}\otimes b)\tilde{T}_{t}^{\ast }
\end{equation}%
for $b\in T_{0}(M)^{\prime }$.
\end{theorem}

\begin{proof}
Suppose $c=\{c_{t}\}_{t\geq 0}$ is in $\mathcal{C}_{q}(E)$, let $H$ be the
Hilbert space of $\{T_{t}\}_{t\geq 0}$, and write
\begin{equation}
\Psi _{t}(b)=\tilde{T}_{t}(c_{t}\otimes b)\tilde{T}_{t}^{\ast }  \label{psit}
\end{equation}%
for $b\in T_{0}(M)^{\prime }$. Then it follows from Proposition~\ref{sub}
that for every $t\geq 0$, $\Psi _{t}$ is a normal, completely positive map
on $T_{0}(M)^{\prime }$ that is subordinate to $\Theta _{t}$. For $t,s\geq 0$
and $b\in T_{0}(M)^{\prime }$, $\Psi _{s}\circ \Psi _{t}(b)=\tilde{T}%
_{s}(c_{s}\otimes (\tilde{T}_{t}(c_{t}\otimes b)\tilde{T}_{t}^{\ast })\tilde{%
T}_{s}^{\ast }=\tilde{T}_{s}(I_{E(s)}\otimes \tilde{T}_{t})(c_{s}\otimes
c_{t}\otimes b)(I_{E(s)}\otimes \tilde{T}_{t}^{\ast })\tilde{T}_{s}^{\ast }=%
\tilde{T}_{s+t}(q_{s+t}\otimes I)(c_{s}\otimes c_{t}\otimes
b)(q_{s+t}\otimes I)\tilde{T}_{s+t}^{\ast }=\tilde{T}_{s+t}(c_{s+t}\otimes b)%
\tilde{T}_{s+t}^{\ast }=\Psi _{s+t}(b)$. (Note that, for every $t$, $\tilde{T%
}_{t}=\tilde{T}_{t}(q_{t}\otimes I_{H})$). Thus $\{\Psi _{t}\}_{t\geq 0}$ is
a semigroup.

Define $S_{t}:E(t)\rightarrow B(H)$ for $t\geq 0$ by setting $S_{t}(\xi
):=T_{t}(c_{t}^{1/2}\xi )$. Then $\{S_{t}\}_{t\geq 0}$ is a covariant
representation of the product system $\{E(t)\}_{t\geq 0}$ on $H$. For $h,k$
in $H$ and a Borel section $t\mapsto \xi _{t}$ of $\{E(t)\}_{t\geq 0}$, the
map $t\mapsto \langle S_{t}(\xi _{t})h,k\rangle =\langle
T_{t}(c_{t}^{1/2}\xi _{t})h,k\rangle $ is Borel because $c$ satisfies (Cq3).
Thus $\{S_{t}\}_{t\geq 0}$ is a measurable representation. It then follows
from Theorem~\ref{cont} that for $b\in T_{0}(M)^{\prime }$, the map $%
t\mapsto \Psi _{t}(b)$ is continuous on $(0,\infty )$. To show continuity at
$t=0$ we need to show that $\cap \{ker(\Psi _{t})|t>0\}=\{0\}$.

For that we first fix $0\leq t\leq s$ and $b\in T_{0}(M)^{\prime }$ and
compute:
\begin{eqnarray}
\Theta _{s-t}\circ \Psi _{t}(b) &=&\tilde{T}_{s-t}(I_{E(s-t)}\otimes (\tilde{%
T}_{t}(c_{t}\otimes b)\tilde{T}_{t}^{\ast }))\tilde{T}_{s-t}^{\ast }
\label{a} \\
&=&\tilde{T}_{s}(q_{s}\otimes I)(I_{E(s-t)}\otimes c_{t}\otimes
b)(q_{s}\otimes I)\tilde{T}_{s}^{\ast }.
\end{eqnarray}%
For $h,k\in H$ and $\xi ,\zeta \in q_{s}E(s)$, $\langle (I_{E(s-t)}\otimes
c_{t}\otimes b)(\xi \otimes h),\zeta \otimes k\rangle =\langle
bh,T_{0}(\langle (I_{E(s-t)}\otimes c_{t})\xi ,\zeta \rangle )k\rangle $
and, using (Cq4), this expression tends to $\langle (I_{E(s)}\otimes b)(\xi
\otimes h),\zeta \otimes k\rangle $ as $t$ tends to $0$. Combining this with
equation (\ref{a}), we see that, for $b\in T_{0}(M)^{\prime }$ and $s>0$,
\begin{equation}
\lim_{t\rightarrow 0}\Theta _{s-t}(\Psi _{t}(b))=\Theta _{s}(b)  \label{b}
\end{equation}%
ultraweakly. Fix $b\in \cap \{ker(\Psi _{t})|t>0\}$ and $s>0$. It follows
from equation (\ref{b}) that $b\in ker(\Theta _{s})$. Since this holds for
all $s>0$, we conclude that $b\in \cap \{ker(\Theta _{t})|t>0\}$. But, from
the continuity of $\Theta _{t}$ it follows that the intersection is just $%
\{0\}$. This proves that $\cap \{ker(\Psi _{t})|t>0\}=\{0\}$ and, using
Theorem~\ref{cont}, we find that $\{\Psi _{t}\}_{t\geq 0}$ is a $CP$%
-semigroup that is subordinate to $\{\Theta _{t}\}_{t\geq 0}$.

For the converse, let $\{\Psi _{t}\}_{t\geq 0}$ be a $CP$-semigroup that is
subordinate to $\{\Theta _{t}\}_{t\geq 0}$. Then by Proposition~\ref{sub}
there is a unique family $c=\{c_{t}\}_{t\geq 0}$ satisfying (Cq1) such that $%
\Psi _{t}$ is of the form (\ref{psit}). It follows from Lemma~\ref{mult}
(iii) (and the uniqueness of the $c_{t}$'s) that, for $s,t\geq 0$, $%
c_{s+t}=q_{s+t}(c_{s}\otimes c_{t})q_{s+t}$. It is left to prove (Cq3) and
(Cq4).

Fix $a,b\in T_{0}(M)^{\prime }$ and $h,k\in H$. Then, for $t\geq 0$, $%
\langle (c_{t}\otimes I_{H})(I_{E(t)}\otimes b)\tilde{T}_{t}^{\ast
}h,(I_{E(t)}\otimes a)\tilde{T}_{t}^{\ast }\rangle =\langle \Psi
_{t}(a^{\ast }b)h,k\rangle $. It results from the continuity properties of $%
\{\Psi _{t}\}_{t\geq 0}$ that this expression is a continuous function of $t$
(on $[0,\infty )$). Since $q_{t}$ is the support projection of $T_{t}$, the
vectors of the form $(I_{E(t)}\otimes b)\tilde{T}_{t}^{\ast }h$, $b\in
T_{0}(M)^{\prime }$ and $h\in H$, span a dense subspace of $q_{t}E(t)\otimes
H$. Thus the map $t\mapsto c_{t}\otimes I_{H}=q_{t}c_{t}q_{t}\otimes I_{H}$
is a Borel section of $\{B(E(t)\otimes H)\}_{t\geq 0}$. Fix Borel sections $%
\{\xi _{t}\}_{t\geq 0}$ and $\{\zeta _{t}\}_{t\geq 0}$ of $\{E(t)\}_{t\geq
0} $. Then for every $h,k\in H$, the map $t\mapsto \langle (c_{t}\otimes
I_{H})(\xi _{t}\otimes h),\zeta _{t}\otimes k\rangle =\langle
h,T_{0}(\langle c_{t}\xi _{t},\zeta _{t}\rangle )k\rangle $ is a Borel map
and (Cq3) follows.

Since both $\{\Theta _{t}\}_{t\geq 0}$ and $\{\Psi _{t}\}_{t\geq 0}$ are $CP$%
-semigroups, it follows from Remark~\ref{weakstrong} (and the fact that
these maps are contractive) that, for every $\omega \in (T_{0}(M)^{\prime
})_{\ast }$ and $s>0$, we have $(\Psi _{t})_{\ast }((\Theta _{s-t})_{\ast
}(\omega ))\rightarrow (\Theta _{s})_{\ast }(\omega )$ in norm as $%
t\rightarrow 0$. Using equation (\ref{a}), we find that, for every $b,a\in
T_{0}(M)^{\prime }$, $s>0$ and $h,k\in H$, $\langle \tilde{T}%
_{s}(q_{s}\otimes I)(I_{E(s-t)}\otimes c_{t}\otimes a^{\ast }b)(q_{s}\otimes
I)\tilde{T}_{s}^{\ast }h,k\rangle \rightarrow \langle \tilde{T}_{s}(I\otimes
a^{\ast }b)\tilde{T}_{s}^{\ast }h,k\rangle $ as $t\rightarrow 0$. Thus $%
\langle (I_{E(s-t)}\otimes c_{t}\otimes I)(I_{E(s)}\otimes b)\tilde{T}%
_{s}^{\ast }h,(I_{E(s)}\otimes a)\tilde{T}_{s}^{\ast }k\rangle \rightarrow
\langle (I_{E(s)}\otimes b)\tilde{T}_{s}^{\ast }h,(I_{E(s)}\otimes a)\tilde{T%
}_{s}^{\ast }k\rangle $. Since vectors of the form $(I_{E(s)}\otimes b)%
\tilde{T}_{s}^{\ast }h$ span a dense subspace of $q_{s}E(s)\otimes H$, it
follows that $q_{s}(I_{E(s-t)}\otimes c_{t})q_{s}\otimes I_{H}\rightarrow
q_{s}\otimes I_{H}$ in the weak operator topology on $B(E(s)\otimes H)$.
Consequently, we have verified (Cq4).
\end{proof}

The referee notes that the algebraic aspects of Theorem \ref%
{subordinatetheorem} may be proved, using duality, from \cite[Theorem 14.3]%
{BS00} or \cite[Theorem A.3]{BBLS04}.

\begin{definition}
\label{cocyle}Let $M$ be a von Neumann algebra and let $\alpha =\{\alpha
_{t}\}_{t\geq 0}$ be an $E$-semigroup of endomorphisms of $M$. A \emph{%
(left) cocycle} for $\alpha $ is is a family of operators $\{C_{t}\}_{t\geq
0}$ in $M$ such that as a function from $[0,\infty )$ to $M$, $%
\{C_{t}\}_{t\geq 0}$ is strongly continuous and such that $%
C_{t+s}=C_{s}\alpha _{s}(C_{t})$ for all $s,t\geq 0$, with $C_{0}=I$. A
cocycle $\{C_{t}\}_{t\geq 0}$ is called \emph{local} if $C_{t}$ commutes
with $\alpha _{t}(M)$ for all $t\geq 0$. We say that $\{C_{t}\}_{t\geq 0}$
is positive and contractive if each of the $C_{t}$ is a positive contraction
operator.
\end{definition}

Of course right cocycles may be defined similarly, but a local left cocycle
is the same as a local right cocycle. The connection between cocycles for
product systems and (local) cocycles for $E$-semigroups is made in the
following proposition.

\begin{proposition}
\label{cocycleequal}Let $\alpha =\{\alpha _{t}\}_{t\geq 0}$ be an $E$%
-semigroup of endomorphisms of a von Neumann algebra acting on a Hilbert
space $H$. Suppose that $\alpha $ is given by an isometric representation $%
\{V_{t}\}_{t\geq 0}$ on $H$ of a measurable product system $\{E(t)\}_{t\geq
0}$ over $M^{\prime }$ through the equation
\begin{equation*}
\alpha _{t}(S)=\tilde{V}_{t}(I_{E(t)}\otimes S)\tilde{V}_{t}^{\ast }\text{,}
\end{equation*}%
$S\in M$. Then for every cocycle for $\{E(t)\}_{t\geq 0}$, $%
c=\{c_{t}\}_{t\geq 0}$ in $\mathcal{C}(E)$, $C=\{C_{t}\}_{t\geq 0}$, defined
by the formula%
\begin{equation}
C_{t}:=\tilde{V}_{t}(c_{t}\otimes I)\tilde{V}_{t}^{\ast }\text{,}
\label{localcocyclequate}
\end{equation}%
$t\geq 0$, is a positive, contractive local cocycle for $\alpha $ on $M$.
Conversely, every positive, contractive local cocycle for $\alpha $ on $M$, $%
\{C_{t}\}_{t\geq 0}$, is given by equation (\ref{localcocyclequate}) for a
unique cocycle $c=\{c_{t}\}_{t\geq 0}$ in $\mathcal{C}(E)$.
\end{proposition}

\begin{proof}
By Lemma \ref{comm}, an element $c=\{c_{t}\}_{t\geq 0}$ in $\mathcal{C}(E)$
gives rise to a local cocycle for $\alpha $, $C=\{C_{t}\}_{t\geq 0}$, via
equation (\ref{localcocyclequate}). The algebraic properties of $C$ are easy
to verify on the basis of this definition. The only possible difficulty is
showing that $C$ is continuous. However, observe that for $S\in M$, $%
C_{t}\alpha _{t}(S)=\tilde{V}_{t}(c_{t}\otimes I)\tilde{V}_{t}^{\ast }\tilde{%
V}_{t}(I_{E(t)}\otimes S)\tilde{V}_{t}^{\ast }=\tilde{V}_{t}(c_{t}\otimes S)%
\tilde{V}_{t}^{\ast }$, i.e., $\{C_{t}\alpha _{t}\}_{t\geq 0}$ is the $CP$%
-semigroup that is subordinate to $\alpha $ determined by $c$ and therefore
is ultraweakly continuous. Since $\{C_{t}\alpha _{t}\}_{t\geq 0}$ is
contractive, ultraweak continuity is the same as continuity in the weak
operator topology since the weak and ultraweak topologies agree on bounded
sets. Thus for each $S\in M$, $t\rightarrow C_{t}\alpha _{t}(S)$ is a weakly
continuous function on $[0,\infty )$ that converges to $S$ weakly as $%
t\rightarrow 0+$. By \cite[Proposition 4.1(1)]{MSQMP}, we conclude that $%
t\rightarrow C_{t}\alpha _{t}(S)$ is strongly continuous on $[0,\infty )$
for each $S\in M$. Thus, in particular for $S=I$, we see that $t\rightarrow
C_{t}\alpha _{t}(I)=\tilde{V}_{t}(c_{t}\otimes I)\tilde{V}_{t}^{\ast }(%
\tilde{V}_{t}\tilde{V}_{t}^{\ast })=C_{t}$ is strongly continuous on $%
[0,\infty )$.

Conversely, if $C=\{C_{t}\}_{t\geq 0}$ is a local cocycle, then $%
\{C_{t}\alpha _{t}\}_{t\geq 0}$ is a $CP$-semigroup that is subordinate to $%
\alpha $. Indeed, the fact that $\{C_{t}\alpha _{t}\}_{t\geq 0}$ is a $CP$%
-semigroup is an easy calculation; the fact that $\{C_{t}\alpha
_{t}\}_{t\geq 0}$ is subordinate to $\alpha $ is simply the observation that
for all $t\geq 0$, $\alpha _{t}(\cdot )-C_{t}\alpha _{t}(\cdot
)=(I-C_{t})\alpha _{t}(\cdot )=(I-C_{t})^{\frac{1}{2}}\alpha _{t}(\cdot
)(I-C_{t})^{\frac{1}{2}}$, which is clearly completely positive. Thus, by
Theorem \ref{subordinatetheorem} there is a unique $c\in \mathcal{C}(E)$
such that $C_{t}\alpha _{t}(S)=\tilde{V}_{t}(c_{t}\otimes S)\tilde{V}%
_{t}^{\ast }$ for all $S\in M$. Letting $S=I$ completes the proof.
\end{proof}

The referee has observed that a purely algebraic form of Proposition \ref%
{cocycleequal} may be deduced from \cite[Lemma 7.5]{BS00} and \cite[Theorem
4.4.8]{BBLS04}.

Note that for a given family of projections $q=\{q_{t}\}_{t\geq 0}$
satisfying the conditions of Definition \ref{reducedcocycles}, the cocycles
for a product system $\{E(t)\}_{t\geq 0}$ in $\mathcal{C}_{q}(E)$, do not
depend upon any particular representation $\{T_{t}\}_{t\geq 0}$. On the
other hand, any representation $\{T_{t}\}_{t\geq 0}$ of $\{E(t)\}_{t\geq 0}$
determines a family $q=\{q_{t}\}_{t\geq 0}$ and then $\mathcal{C}_{q}(E)$
depends only on the product system and the spaces $\{Ker(T_{t})\}$.
Consequently, we have the following

\begin{corollary}
\label{orderisom}Let $\{E(t)\}_{t\geq 0}$ be a measurable product system and
let $\{T_{t}\}_{t\geq 0}$ and $\{S_{t}\}_{t\geq 0}$ be two measurable
representations such that for each $t\geq 0$ the support projection of $%
T_{t} $, $q_{T_{t}}$, equals the support projection of $S_{t}$,$\ q_{S_{t}}$%
, and such that the semigroups $\Theta ^{T}$ and $\Theta ^{S}$ associated
with $\{T_{t}\}_{t\geq 0}$ and $\{S_{t}\}_{t\geq 0}$, as in equation (\ref%
{thethat}), are $CP$-semigroups. Then there is an order isomorphism between
the set consisting of those $CP$-semigroups that are subordinate to $\Theta
^{T}$ and the set consisting of those that are subordinate to $\Theta ^{S}$.
\end{corollary}

Using the analysis of \cite{MSQMP} and the uniqueness of the minimal
endomorphic dilation (\cite[Section 8.9 ff.]{Arv03}), we can now prove

\begin{corollary}
\label{dilation}Let $\Theta $ be $CP_{0}$-semigroup on a countably
decomposable von Neumann algebra $N$, acting on the Hilbert space $H$, let $%
\alpha $ be its minimal endomorphic dilation acting on the von Neumann
algebra $R$, which is contained in $B(K)$, and let $W$ be the embedding of $%
H $ into $K$. Then there is an order isomorphism between the $CP$-semigroups
subordinate to $\Theta $ and those that are subordinate to $\alpha $.
Further, if $\Psi $ is subordinate to $\alpha $, then its compression to $H$
is subordinate to $\Theta $.
\end{corollary}

\begin{proof}
Associated with $\Theta $, we get a product system $E_{\Theta }$ as
described in equation (\ref{dirlim}) and the discussion preceding it. We
also get a fully coisometric representation $\{T_{t}\}_{t\geq 0}$ of $%
E_{\Theta }$ (called the identity representation in \cite[Theorem 3.9]{MSQMP}%
). Using Theorem \ref{dilrep}, we dilate the identity representation to a
representation $\{V_{t}\}_{t\geq 0}$ that is both isometric and fully
coisometric. The minimal endomorphic dilation of $\Theta $ is then defined
by setting $\alpha _{t}(b)=\tilde{V}_{t}(I_{E(t)}\otimes b)\tilde{V}%
_{t}^{\ast }$ for $b\in V_{0}(M)^{\prime }$ (where $M=T_{0}(N)^{\prime }$).
Since $\{V_{t}\}_{t\geq 0}$ is isometric and fully coisometric, the order
interval of subordinates of $\{\alpha _{t}\}_{t\geq 0}$ is order isomorphic
to $\mathcal{C}(E_{\Theta })$ by Theorem ~\ref{subordinatetheorem}. Also by
this theorem, the order interval of subordinates of $\Theta $ is order
isomorphic to $\mathcal{C}_{q}(E_{\Theta })$, where $q=\{q_{t}\}_{t\geq 0}$
is the family of support projections of $\{T_{t}\}_{t\geq 0}$. So, to
complete the proof, we have to establish an order isomorphism $\Gamma $ from
$\mathcal{C}(E_{\Theta })$ onto $\mathcal{C}_{q}(E_{\Theta })$. We simply
define
\begin{equation*}
\Gamma (c)=qcq
\end{equation*}%
where $qcq:=\{q_{t}c_{t}q_{t}\}_{t\geq 0}$. It is easy to check that $\Gamma
$ maps $\mathcal{C}(E_{\Theta })$ into $\mathcal{C}_{q}(E_{\Theta })$ and
preserves the order. What we need to prove is that it is injective and
surjective; that is, we need to show that given $c\in \mathcal{C}%
_{q}(E_{\Theta })$, there is a unique $c^{\prime }\in \mathcal{C}(E_{\Theta
})$ such that $qc^{\prime }q=c$. So fix $c\in \mathcal{C}_{q}(E_{\Theta })$.
For every partition $\mathcal{P}=\{0=t_{0}<\cdots <t_{n}=t\}$ of $[0,t]$ we
define
\begin{equation*}
c_{\mathcal{P},t}=c_{t_{1}}\otimes c_{t_{2}-t_{1}}\otimes \cdots \otimes
c_{t-t_{n-1}}
\end{equation*}%
and, similarly, we define
\begin{equation*}
q_{\mathcal{P},t}=q_{t_{1}}\otimes q_{t_{2}-t_{1}}\otimes \cdots \otimes
q_{t-t_{n-1}}.
\end{equation*}%
It follows from (Cq1) and (Cq2) that
\begin{equation*}
c_{\mathcal{P},t}=q_{\mathcal{P},t}c_{\mathcal{P},t}q_{\mathcal{P},t}\text{.}
\end{equation*}%
It also follows from (Cq1) and (Cq2) that if $\mathcal{P}^{\prime }$ refines
$\mathcal{P}$, then
\begin{equation}
c_{\mathcal{P},t}=q_{\mathcal{P},t}c_{\mathcal{P}^{\prime },t}q_{\mathcal{P}%
,t}\text{.}  \label{refine}
\end{equation}%
Comparing \cite[Lemma 4.3(2)]{MSQMP} with the definition of $q_{t}$ as the
support projection of $T_{t}$, we see that $q_{t}$ is the projection from $%
E_{\Theta }(t)$ onto $\mathcal{E}_{t}$. Thus, for every partition $\mathcal{P%
}$ of $[0,t]$, $q_{\mathcal{P},t}$ is the projection onto $\mathcal{E}(%
\mathcal{P},t)$ and it follows from equation (\ref{dirlim}) that $q_{%
\mathcal{P},t}$ converges ($\sigma -$weakly) to $I_{E_{\Theta }(t)}$. It now
follows from equation (\ref{refine}) that there is a unique positive
contraction $c_{t}^{\prime }\in \mathcal{L}(E_{\Theta }(t))\cap \varphi
_{t}(M)^{\prime }$ that satisfies
\begin{equation*}
c_{\mathcal{P},t}=q_{\mathcal{P},t}c_{t}^{\prime }q_{\mathcal{P},t}.
\end{equation*}%
for all partitions. In particular, $q_{t}c_{t}^{\prime }q_{t}=c_{t}$; that
is $\Gamma (c^{\prime })=c$. It is straightforward to check that $c^{\prime
} $ lies in $\mathcal{C}(E_{\Theta })$. If $d\in \mathcal{C}(E_{\Theta })$
also satisfies $q_{t}d_{t}q_{t}=c_{t}$, then $q_{\mathcal{P},t}d_{t}q_{%
\mathcal{P},t}=c_{\mathcal{P},t}=q_{\mathcal{P},t}c_{t}^{\prime }q_{\mathcal{%
P},t}$ for all partitions of $[0,t]$. Taking the limit over all partitions,
we find that $d_{t}=c_{t}^{\prime }$, completing the proof of the first
assertion.

For the second, simply note that if $\Psi $ is subordinate to $\alpha $, so
that it is given by the formula
\begin{equation*}
\Psi _{t}(b)=\tilde{V}_{t}(c_{t}\otimes b)\tilde{V}_{t}^{\ast }\text{,}
\end{equation*}%
$b\in R$, then the formula
\begin{multline*}
a\rightarrow W^{\ast }\Psi _{t}(WaW^{\ast })W=W^{\ast }\tilde{V}%
_{t}(I\otimes W)(c_{t}\otimes a)(I\otimes W^{\ast })\tilde{V}_{t}^{\ast }W \\
=\tilde{T}_{t}(c_{t}\otimes a)\tilde{T}_{t}^{\ast }=\tilde{T}%
_{t}(q_{t}c_{t}q_{t}\otimes a)\tilde{T}_{t}^{\ast }\text{,}
\end{multline*}%
$a\in N$, defines a $CP$-semigroup that is subordinate to $\Theta $.
\end{proof}

As we mentioned in the introduction, this result improves upon a theorem of
Bhat an Skeide \cite[Theorem 14.3]{BS00} and on theorems of Bhat \cite[%
Section 5]{bB01} and Powers \cite[Theorems 3.4 and 3.5]{Po03}.

Finally, it is interesting to note the following result, which is related to
\cite[Theorem 4.3]{bB02} and which sharpens Theorem \ref{subsys}.

\begin{theorem}
\label{subpro}Let $N$ be a von Neumann algebra represented faithfully on a
Hilbert space $H$, and let $\alpha =\{\alpha _{t}\}_{t\geq 0}$ be an $E$%
-semigroup acting on $N$. Let $E=\{E(t)\}_{t\geq 0}$ be the product system
for $\alpha $ and let $\Theta =\{\Theta _{t}\}_{t\geq 0}$ be a $CP$%
-semigroup that is subordinate to $\alpha $. If $c\in \mathcal{C}(E)$ is the
cocycle guaranteed by Theorem \ref{subordinatetheorem} such that $\Theta
=\alpha _{c}$ and if $q_{t}$ is the projection onto the closure of the range
of $c_{t}^{1/2}$, then $q_{t+s}=q_{t}\otimes q_{s}$ for all $t,s\geq 0$ and
the product system $\{E_{\Theta }(t)\}_{t\geq 0}$ is isomorphic to the
product system $\{q_{t}E(t)\}_{t\geq 0}$.
\end{theorem}

\begin{proof}
Let $V=\{V_{t}\}_{t\geq 0}$ be the isometric representation of $%
\{E(t)\}_{t\geq 0}$ that gives rise to $\{\alpha _{t}\}_{t\geq 0}$. Applying
Theorem~\ref{subordinatetheorem}, we see that $c,$ $V$ and $\Theta $ are
related by the equation
\begin{equation*}
\Theta _{t}(a)=\tilde{V}_{t}(c_{t}\otimes a)\tilde{V}_{t}^{\ast }\text{,}
\end{equation*}%
which is valid for all $a\in N$ and $t\geq 0$. (Note that since each $V_{t}$
is injective, $c$ is assumed to be a cocycle of the product system and not a
reduced cocycle.) Hence $\{\Theta _{t}\}_{t\geq 0}$ is the semigroup
associated with the covariant representation $\{T_{t}\}_{t\geq 0}$ where $%
T_{t}:=V_{t}\circ c_{t}^{1/2}$. Using equation (\ref{q}), we find that the
support projection of $T_{t}$, $q_{T_{t}}$, satisfies the equation
\begin{equation}
q_{T_{t}}E(t)\otimes H=[(I_{E(t)}\otimes N)(c_{t}^{1/2}\otimes I)\tilde{V}%
_{t}^{\ast }H]=\overline{(c_{t}^{1/2}\otimes I)(E(t)\otimes H)}.
\label{Range}
\end{equation}%
Thus $q_{T_{t}}$ is the projection onto $\overline{c_{t}^{1/2}E(t)}$, i.e., $%
q_{T_{t}}=q_{t}$. Since $c_{t+s}=c_{t}\otimes c_{s}$, we have $%
q_{t+s}=q_{t}\otimes q_{s}$ and we may apply Theorem~\ref{subsys} to
complete the proof.
\end{proof}

Theorem \ref{subpro} has the following immediate corollary that we promised
in Remark \ref{endprod}.

\begin{corollary}
\label{PersistentProductSyst}If $\Theta =\{\Theta _{t}\}_{t\geq 0}$ is a $CP$%
-semigroup on a von Neumann algebra $N$ that is subordinate to an $E$%
-semigroup, then the Arveson-Stinespring correspondences $\{\mathcal{E}%
_{t}\}_{t\geq 0}$ associated with $\Theta $ form a multiplicative family,
i.e., $\mathcal{E}_{t+s}\simeq \mathcal{E}_{t}\otimes \mathcal{E}_{s}$, for
all $s,t\geq 0$.
\end{corollary}


\begin{thebibliography}{99}
\bibitem{AnD90} C. Anantharaman-Delaroche, \emph{On completely positive maps
defined by an irreducible correspondence}, Canad. Math. Bull. 33 (1990),
434-441.

\bibitem{Arv69} W.B. Arveson, \emph{Subalgebras of }$C^{\ast }$\emph{%
-algebras}, Acta Math. 123 (1969), 141-224.

\bibitem{Arv89} W.B. Arveson, \emph{Continuous analogues of Fock space},
Mem. Amer. Math. Soc. \textbf{80} (1989).

\bibitem{Arv97} W. B. Arveson, \emph{The index of a quantum dynamical
semigroup, }J. Functional Anal. \textbf{146} (1997), 557-588.

\bibitem{Arv03} W.B. Arveson, \textit{Noncommutative Dynamics and
E-Semigroups}, Springer Monographs in Mathematics, Springer-Verlag New York
Berlin Heidelberg, 2003.

\bibitem{BDH88} M. Baillet, Y. Denizeau and J.-F. Havet, \emph{Indice d'une
esperance conditionelle}, Comp. Math. \textbf{66} (1988), 199-236.

\bibitem{BBLS04} S. Barreto, B. V. R. Bhat, V. Liebscher and M. Skeide,
\emph{Type I product systems of Hilbert modules}, J. Funct. Anal. \textbf{212%
} (2004), 121--181.

\bibitem{bB01} B. V. R. Bhat, \emph{Cocycles of CCR flows}, Mem. Amer. Math.
Soc. \textbf{149} (2001), no. 709, x+114 pp.

\bibitem{bB02} B. V. R. Bhat, \emph{Minimal isometric dilations of operator
cocycles}, Integral Equations Operator Theory \textbf{42} (2002), 125--141.

\bibitem{BS00} B. V. R. Bhat and M. Skiede, \emph{Tensor product systems of
Hilbert modules and dilations of completely positive semigroups}, Inf. Dim.
Anal. Quantum Prob. and Rel. Topics \textbf{3} (2000), 519-575.

\bibitem{BGR77} L. Brown, P. Green and M. Rieffel, \emph{Stable isomorphism
and strong Morita equivalence of} $C^{\ast }$-\emph{algebras},.Pacific J.
Math. \textbf{71} (1977), 349--363.

\bibitem{Dix} J. Dixmier, \textit{Von Neumann algebras}, North-Holland
Mathematical Library, North-Holland, 1981.

\bibitem{E65} E.G. Effros, \emph{The Borel structure of von Neumann algebras
on a separable Hilbert space}, Pac. J. Math. \textbf{15} (1965), 1153-1164.

\bibitem{EK98} D. Evans and Y. Kawahigashi, \emph{Quantum symmetries on
operator algebras}, Oxford Mathematical Monographs. Oxford Science
Publications. The Clarendon Press, Oxford University Press, New York, 1998.
xvi+829 pp.

\bibitem{FD88} J. M. G. Fell and R. S. Doran, \emph{Representations of $\ast
$-algebras, locally compact groups, and Banach $\ast $-algebraic bundles.}
Pure and Applied Mathematics, 125 and 126. Academic Press, Inc., Boston, MA,
1988.

\bibitem{nF02} N. Fowler, \emph{Discrete product systems of Hilbert bimodules%
}, Pacific J. Math. \textbf{204} (2002), 335--375.

\bibitem{GLR85} P. Ghez, R. Lima, and J. Roberts, $W^{\ast }$-\emph{%
categories}, Pacific J. Math. \textbf{120} (1985), 79--109.

\bibitem{HP57} E. Hille and R. Phillips, \textit{Functional Analysis and
Semi-Groups}, Third printing of the revised edition of 1957, American
Mathematical Society Colloquium Publications, Vol. XXXI. Amer. Math. Soc.,
Providence, RI, 1974.

\bibitem{iH04} I. Hirshberg, \emph{$C^{\ast }$-Algebras of Hilbert module
product systems}, J. Reine Angew. Math. \textbf{570} (2004), 131-142.

\bibitem{KS97} M. Khoshkam and G. Skandalis, \emph{Toeplitz algebras
associated with endomorphisms and Pimsner-Voiculescu exact sequences,}
Pacific J. Math. \textbf{181} (1997), 315--331.

\bibitem{L94} E.C. Lance , \emph{Hilbert $C^{\ast }$-modules, A toolkit for
operator algebraists}, London Math. Soc. Lecture Notes series \textbf{210}
(1995). Cambridge Univ. Press.

\bibitem{M03} D. Markiewicz, \emph{On the product system of a completely
positive semigroup}, J. Funct. Anal. \textbf{200} (2003), 237-280.

\bibitem{Mi89} J. Mingo, \emph{The correspondence associated to an inner
completely positive map}, Math. Ann. \textbf{284} (1989), 121-135.

\bibitem{pM01} P.S. Muhly, \emph{Bundles over groupoids }in Groupoids in
analysis, geometry, and physics (Boulder, CO, 1999), 67--82, Contemp. Math.,
\textbf{282}, Amer. Math. Soc., Providence, RI, 2001.

\bibitem{MSTensor} P.S. Muhly and B. Solel, \emph{Tensor algebras over }$%
C^{\ast }$\emph{-correspondences : Representations, dilations and }$C^{\ast
} $\emph{-envelopes}, J. Funct. Anal. \textbf{158} (1998), 389-457.

\bibitem{MSQMP} P.S. Muhly and B. Solel, \emph{Quantum Markov processes
(correspondences and dilations)}, Int. J. of Math. 13(2002), 863-906.

\bibitem{MSCurv} P.S. Muhly and B. Solel, \emph{The curvature and index of
completely positive maps}, Proc. London Math. Soc. 87 (2003), 748-778.

\bibitem{MSHardy} P.S. Muhly and B. Solel, \emph{Hardy algebras, }$W^{\ast }$%
\emph{-correspondences and interpolation theory}, Math. Annalen \textbf{330}
(2004), 353--415.

\bibitem{Pa73} W.L. Paschke, \emph{Inner product modules over }$B^{\ast }$%
\emph{-algebras}, Trans. Amer. Math. Soc. 182 (1973), 443-468.

\bibitem{P86} S. Popa, \emph{Correspondences}, Preprint (1986).

\bibitem{rP88} R. T. Powers, \emph{An index theory for semigroups of }$\ast $%
\emph{-endomorphisms of }$\mathfrak{B}(\mathfrak{H})$\emph{\ and type }$%
II_{1}$\emph{\ factors}, Can. J. Math. \textbf{40} (1988), 86--114.

\bibitem{Po03} R.T. Powers, \emph{Continuous spatial semigroups of
completely positive maps of }$\mathfrak{B}(\mathfrak{H})$, New York J. Math.
9 (2003), 165-269.

\bibitem{R74} M.A. Rieffel, \emph{Induced representations of $C^{\ast }$%
-algebras}, Adv. in Math. 13 (1974), 176-257.

\bibitem{mS03} M. Skeide, \emph{Commutants of von Neumann modules,
representations of }$\mathfrak{B}^{a}(E)$\emph{\ and other topics related to
product systems of Hilbert modules}, Advances in quantum dynamics (South
Hadley, MA, 2002), 253--262, Contemp. Math., \textbf{335}, Amer. Math. Soc.,
Providence, RI, 2003.

\bibitem{mS03a} M. Skeide, \emph{Dilation theory and continuous tensor
product systems of Hilbert modules}, in \emph{Quantum probability and
infinite dimensional analysis (Burg, 2001)}, 215--242, QP--PQ: Quantum
Probab. White Noise Anal., 15, World Sci. Publishing, River Edge, NJ, 2003.

\bibitem{mSp03} M. Skeide, \emph{von Neumann Modules, Intertwiners and
Self-Duality}, J. Operator Theory \textbf{54} (2005), 119--124.

\bibitem{mS06} M. Skeide, \emph{Isometric dilations of representations of
product systems via commutants}, Preprint, ArXiv:math.OA/0602459.
\end{thebibliography}
\end{document}